\documentclass[10pt,a4paper,twoside]{amsart}

\usepackage{amssymb}
\usepackage{graphicx}
\usepackage{bm}
\usepackage[usenames,dvipsnames]{xcolor}
\usepackage{cite}

\newcommand{\rmd}{\mathrm{d}}

\newcommand{\mFr}{\mathfrak{m}}
\newcommand{\dFr}{\mathfrak{d}}

\newcommand{\sfit}[1]{\textsf{\textit{#1}}}

\newcommand{\Jac}{\mathrm{Jac}}

\DeclareMathOperator{\Complex}{\mathbb{C}}
\DeclareMathOperator{\Integer}{\mathbb{Z}}
\DeclareMathOperator{\Natural}{\mathbb{N}}
\DeclareMathOperator{\wgt}{wgt}

\DeclareMathOperator{\rank}{rank}

\DeclareMathOperator{\Resultant}{Resultant}
\DeclareMathOperator{\inv}{Inv}
\DeclareMathOperator{\add}{Add}

\theoremstyle{plain}
\newtheorem{lem}{Lemma}
\newtheorem{theo}{Theorem}
\newtheorem{prop}{Proposition}
\newtheorem*{prop*}{Proposition}

\theoremstyle{definition} 
\newtheorem{defn}{Definition}
\newtheorem{rem}{Remark}
\newtheorem{exam}{Example}

\title[]{Addition via reduction algorithm on~trigonal~curves}
\author{J Bernatska, Y Kopeliovich}
\address{}
\email{}
\date{\today}

\allowdisplaybreaks[4]

\begin{document}
 
\maketitle 
\begin{abstract}
In this paper we propose a direct and explicit realization of addition of divisors by means of 
an iterative reduction algorithm. Each iteration of the algorithm is the
reduction of a degree $g+1$ divisor to a divisor of degree~$g$.
Such an approach allows to carry out all computations explicitly in a symbolic form,
which is done for curves $C_{3,4}$, $C_{3,5}$ in this paper, and also for curves of higher genera
up to $C_{3,13}$, $C_{3,14}$.
\end{abstract}

\section{Introduction}
Cryptosystems based on the discrete-log problem on Jacobian groups continue to draw attention.
The focus has shifted from elliptic and hyperelliptic curves towards
superelliptic and more general algebraic curves.
As stated in \cite{Bas2},  higher genera ($>4$) are considered 
as less secure than elliptic for the same group order, and curves of genus $3$ 
are the most attractive ones from the cryptographic point of view.
At the  same time, increasing the genus of a curve allows to 
decrease the size of the ground field for the same order of magnitude.
Our results are not restricted to genera $3$ and $4$ only, but  applicable to a trigonal curve of 
an arbitrary genus.

We start with several definitions.
\begin{defn}
A divisor of degree $0$ is called a \emph{reduced divisor} on a genus $g$ algebraic curve if 
it has the form: $\widetilde{D} - p \infty$ and $\deg \widetilde{D} = p$, $p \leqslant g$.
Such a divisor is non-special if $p=g$, and special if $p<g$.
\end{defn} 

It follows from the Riemann-Roch theorem that any non-special divisor of the form $D-(\deg D)\infty$
is equivalent to a reduced non-special divisor. 
An accurate explanation of what we call a non-special divisor is given in subsection~\ref{ss:NSpDiv}.
Below we deal with  reduced non-special divisors only, and call them simply reduced divisors.

\newtheorem*{Rproblem}{Reduction Problem}
\begin{Rproblem}
Given a non-special divisor $D - (\deg D)\infty$ of degree $g+m$ with $m>0$ 
on an algebraic curve of genus $g$, find the corresponding reduced divisor $\widetilde{D}$ 
such that $D - {(g+m) \infty}$ is equivalent to $\widetilde{D} - g \infty$.
\end{Rproblem} 
\newtheorem*{Aproblem}{Addition Problem}
\begin{Aproblem} 
Given two non-special divisors $D_1$ and 
$D_2$ of degrees $g+m_1$ and $g+m_2$, $m_1,\,m_2\,{\geqslant}\, 0$, respectively,
find a reduced divisor $\widetilde{D}$ such that $D_1+D_2 - (2g+m_1+m_2)\infty \sim \widetilde{D}-g\infty$. 
\end{Aproblem}

Evidently, solving the reduction problem will solve the addition problem. 
Indeed, we can assume that $D_1$, $D_2$ together compose a non-special divisor $D$ of degree $2g+m_1+m_2$, 
so we come to the reduction problem for the new divisor $D$.
On the other hand, the standard addition problem arises when the both divisors $D_1$, $D_2$ are 
firstly reduced to divisors $\widetilde{D}_1$, $\widetilde{D}_2$ of degree $g$ each.

Like many authors, we employ the class of $C_{a,b}$ curves, which have
one cusp isomorphic to the origin of an affine plane curve $x^a+y^b$ as a singular point. 
This type of curves contains the maximal number of rational places 
on a curve \cite{Miu93}, and became popular in cryptographic literature: 
\cite{KM, Ari99, Ari03, HaSu00, KM1,KM2, Bas1}.
The notion of  $C_{a,b}$ curve coincides with the $(n,s)$-curve introduced in \cite{bel99}. The latter approach 
is more preferable for us due to the significant results in the theory of addition laws 
on Jacobian varieties of such curves presented in \cite{bl05}.

The history of cryptosystems on non-elliptic algebraic curves started from \cite{Kob2,Ca}.
Many improvements of the fast addition algorithm were suggested
mostly in genera $2$ and $3$, for example  \cite{Lan, MCT, KGMCT}. 
Less is done on fast addition in non-hyperelliptic cases. 
Some papers and conference talks consider special classes of non-hyperelliptic curves of genus~$3$. 
These are Picard curves, also called superelliptic cubics, 
which are cyclic trigonal curves of genus~$3$:  \cite{Bar1,Bar2, FlonOyo, Bas2}, 
and a more general trigonal curve $C_{3,4}$: \cite{FlonOyo2, Bas1,Ari03,KM}. 
All the papers realize the same algorithm consisting of two steps: (i) addition of two non-special divisors, 
(ii) reduction to the equivalent reduced divisor. 

One of the goals of the present paper is to explain addition on
Jacobians of algebraic curves in a clear and simple language,
and make it possible for a non-specialist to understand and implement it.

Our approach is based on the theory of addition laws from \cite{bl05},
which gives an easy receipe how to construct functions defining a divisor on a curve,
how many functions are needed, and which weight (the number of zeros) they should have.
In particular, it is known that two functions of weights $2g$ and $2g+1$ on a curve of genus $g$
define a non-special divisor of degree $g$ unambiguously.

Another problem which arises in this connection is that
the two functions which define a divisor can not be chosen arbitrary.
They have to be constructed so that their divisors of zeros on a curve intersect 
in the required divisor. 
None of the papers addresses such a problem, and none gives a receipe how to construct such functions. 

Finally, the proposed reduction algorithm is a universal solution of the
problem of addition on Jacobians. We suggest to add  divisors point by point. 
This makes the algorithm independent of the degrees of divisors.
Therefore, all computations can be done  explicitly, and before a machine realisation.

One of possible applications of the algorithm can be in finding conditions on torsion points that exist in Jacobians similar to the
mentioned in \cite{Ca2} and \cite{Uc1,Uc2}. Those conditions were used in \cite{DeJMu} to give new algorithms to compute Falting's invariants for hyperelliptic curves. 
Similar results should be available if a corresponding theory of addition is developed for an arbitrary curve.

This paper is organized as follows. In Preliminaries we give 
a brief explanation which types of curves we use,
and recall the recent results in the theory of addition laws
on algebraic curves from \cite{bl05}. 
Section~\ref{s:RatF} is devoted to entire rational functions on an  algebraic curve,
which are used to define divisors. We work with effective divisors only, and so
the class of entire rational functions fullfills our needs completely. 
In Section~\ref{s:Main} all steps of the proposed 
reduction algorithm are described in detail. 
The algorithm itself is presented in Section~\ref{s:RedAlg}.
Finally,  we give an example of implementation of the algorithm on a genus $3$ trigonal curve
in Section~\ref{s:C34impl}. Another example on a genus $4$ trigonal curve
is given in Appendix.

\section{Preliminaries}
\subsection{$(n,s)$-Curves}\label{ss:nsC}
An \emph{$(n,s)$-curve},  introduced in  \cite{bel99}, is defined by the equation
\begin{subequations}\label{CnsEq}
\begin{align}
0=&f(x,y) = y^n + x^s + \sum_{j=0}^{n-2} \sum_{i=0}^{s-2} \lambda_{ns - sj-ni}  y^j x^i,\\
&\lambda_{k\leqslant 0} = 0,  \label{lambdaCond}
\end{align}
\end{subequations}
with co-prime integers $n$ and $s$, and $(x,y)\in \Complex^2$, $\lambda_k \in \Complex$.
We use notation $C_{n,s}$ for such a curve.
Function $f$ is the universal unfolding of the Pham singularity $y^n + x^s=0$ 
with the minimal number of parameters $\lambda_k \in \Complex$.
This type of curves serves as a kind of canonical forms among plane algebraic curves, 
like Weierstrass form serves as the canonical form of elliptic curves. In other words,
any plane algebraic curve is mapped by a M\"{o}bius transformation into an $(n,s)$-curve, 
possibly with double points.
As seen from \eqref{CnsEq}, 
an $(n,s)$-curve is constructed in such a way that the infinity is 
a Weierstrass point, and
a branch point connecting all $n$ sheets of the curve. The infinity serves as the base point in Abel's map.

As shown in  \cite{bel99}, the genus of an $(n,s)$-curve is computed by the formula 
\begin{equation}\label{Genus}
g = \frac{1}{2} (n-1)(s-1),
\end{equation}
which is provided by condition  \eqref{lambdaCond}. 
The curve \eqref{CnsEq} is supposed to be non-degenerate, that is its genus equals $g$. 
In the case of a genus less $g$, the corresponding non-degenerate $(n,s)$-form
should be used.

We also introduce the notion of the \emph{S\={a}to weight},
which is respected by the theory of $(n,s)$-curves. 
The S\={a}to weight equals the opposite to the exponent of the leading term in the expansion near infinity. 
Actually, parametrisation of  \eqref{CnsEq} is 
\begin{equation}\label{param}
x=\xi^{-n},\qquad y = \xi^{-s}(1+O(\lambda)),
\end{equation}
where $\xi$ serves as a local parameter  near infinity, and $\lambda$ denotes 
the collection of all parameters of the $(n,s)$-curve in question. 
Thus, the S\={a}to weights of $x$ and $y$ are
 $\wgt x = n$, and $\wgt y = s$. The weight is also assigned to every function, 
 for example, $\wgt f = n s$.
The weights are used as indices of parameters $\lambda_k$, namely: $\wgt \lambda_k = k$.  
Note, that $\lambda_k$ in the equation of an $(n,s)$-curve have only positive weights,
parameters $\lambda_{k\leqslant 0}$ of non-positive weights are supposed to be zero.
With the help of the S\={a}to weights an order relation is introduced in the space of monomials $y^j x^i$.
The latter are used for constructing the equation of a curve 
and also entire rational functions on it. 

The Weierstrass gap sequence $\mathfrak{w}$ 
is obtained as the complement to the following 
sequence in the set of natural numbers:
\begin{gather*}
 \mathfrak{w}^\ast = a n + b s, \quad a,b = \{0\} \cup \Natural.
\end{gather*}
 Thus $\mathfrak{w} = \Natural \backslash \mathfrak{w}^\ast$.

\subsection{Trigonal curves}
We use the theory of $(n,s)$-curves to classify trigonal curves, see \cite{bel00}.
There are two canonical types of trigonal curves: 
$(3,3\mFr+1)$-curves of genera $3\mFr$, denoted below by $\mathfrak{C}^1$, 
and $(3,3\mFr+2)$-curves of genera $3\mFr+1$, denoted by $\mathfrak{C}^2$.
Let a trigonal curve be defined by the equation 
\begin{subequations}\label{C3def}
\begin{gather}\label{C3Eeq}
0=f(x,y) = y^3 + y^2 \mathcal{T}(x) + y \mathcal{Q}(x) + \mathcal{P}(x), 
\end{gather}
where
\begin{align}
&\mathfrak{C}^1:& &\mathcal{P}(x) = x^{3\mFr+1} + \sum_{i=0}^{3\mFr} \lambda_{9\mFr+3-3i} x^i,& \\
&&  &\mathcal{Q}(x) = \sum_{i=0}^{2\mFr} \lambda_{6\mFr+2-3i} x^i, \qquad 
\mathcal{T}(x) =  \sum_{i=0}^{\mFr} \lambda_{3\mFr+1-3i} x^i;& \notag \\
&\mathfrak{C}^2:& &\mathcal{P}(x) = x^{3\mFr+2} + \sum_{i=0}^{3\mFr+1} \lambda_{9\mFr+6-3i} x^i,& \\
&& &\mathcal{Q}(x) = \sum_{i=0}^{2\mFr+1} \lambda_{6\mFr+4-3i} x^i,\qquad 
\mathcal{T}(x) =  \sum_{i=0}^{\mFr} \lambda_{3\mFr+2-3i} x^i.& \notag
\end{align}
\end{subequations}
All curves are supposed to be non-degenerate. Here we extended an $(n,s)$-curve \eqref{CnsEq}
by the extra terms: $x^g$ and $y^2 x^i$ with $i=0$, \ldots $\mFr$, 
which increase the complexity of a curve in the view of cryptographic applications.

The simplest trigonal curve of type $\mathfrak{C}^1$ is  $C_{3,4}$ of genus $g=3$, namely
\begin{multline}\label{C34eq1}
 0 = f(x,y)  = y^3 + y^2 (\lambda_1 x + \lambda_4 )+ y (\lambda_2 x^2 + \lambda_5 x + \lambda_8) \\
  + x^4 +\lambda_3 x^3 + \lambda_6 x^2 + \lambda_9 x + \lambda_{12},
\end{multline}
and  of type $\mathfrak{C}^2$ is $C_{3,5}$ of genus $g=4$  
\begin{multline}\label{C35eq1}
 0 = f(x,y)  = y^3  + y^2 (\lambda_2 x + \lambda_5) + y (\lambda_1 x^3 + \lambda_4 x^2 + \lambda_7 x + \lambda_{10}) \\
  + x^5 + \lambda_3 x^4 +\lambda_6 x^3 + \lambda_9 x^2 + \lambda_{12} x + \lambda_{15} .
\end{multline}

 We distinguish between $\mathfrak{C}^1$ and $\mathfrak{C}^2$  
 because of the differences in  gap sequences and 
 in the order of monomials. 
The gap sequences (supposed to be sorted in the ascending order with respect to the S\={a}to weight) are
\begin{align}\label{GSeqTrig}
\begin{split}
\mathfrak{C}^1:&\quad\mathfrak{w}=\{3k-2 \mid k=1,\dots, \mFr\} \cup \{3k-1 \mid k=1,\dots, 2\mFr\}; \\
\mathfrak{C}^2:&\quad\mathfrak{w}=\{3k-1 \mid k=1,\dots, \mFr\} \cup \{3k-2 \mid k=1,\dots, 2\mFr+1\}.
\end{split}
\end{align}
The ordered lists of monomials are
\begin{subequations}\label{Monoms}
\begin{align}
&\mathfrak{C}^1:\quad \mathfrak{M}    \label{Monom1}
= \big\{1,  x, \dots, x^{\mFr-1}, x^\mFr, y, x^{\mFr+1}, y x, \dots, x^{2\mFr-1}, y x^{\mFr-1},\\ 
&\qquad x^{2\mFr},  y x^\mFr,  y^2,  x^{2\mFr+1},  y x^{\mFr+1}, 
\{ y^2 x^i, x^{2\mFr+1+i},  y x^{\mFr+1+i} \mid i\in\Natural \} \big\}, \notag  \\  
&\mathfrak{C}^2:\quad \mathfrak{M} \label{Monom2}
= \big\{ 1,  x,  \dots, x^{\mFr-1},  x^\mFr,  y,  x^{\mFr+1},  
y x,  \dots, x^{2\mFr-1}, y x^{\mFr-1}, x^{2\mFr}, \\ 
&\qquad y x^\mFr, x^{2\mFr+1},  y^2,  yx^{\mFr+1}, x^{2\mFr+2},  
\{y^2 x^i,  yx^{\mFr+1+i}, x^{2\mFr+2+i} \mid i\in\Natural \} \big\}.  \notag
\end{align}
\end{subequations}

 In general, such a list  serves as a characteristic of an $(n,s)$-curve.
The first $g$ monomials appear as numerators of holomorphic differentials.

\begin{exam}\label{E:CEq}
Below we perform numerical computations on the curve
\begin{equation}\label{C34Ex1}
y^3  + y^2 (2 x-1) - y(4x^2 + 3x + 2) + x^4 - 49 x^3 +197 x^2 - 52 x - 334 = 0,
\end{equation}
specially designed to contain at least four integer points: $(-1,\,5)$, $(3,\,1)$, $(4,-3)$, and $(2,\,-1)$.
We use these points to simulate a solution of the standard cryptography problem
 within $\Integer$.
\end{exam}

\subsection{Jacobian variety and Abel's map}
Each plane algebraic curve  $C$ is related to a Jacobian variety 
$\Jac (C)$ with coordinates $u=(u_{w_1},u_{w_2},\dots, u_{w_g})$,
where $\{w_k\}_{k=1}^g$ is the gap sequence~$\mathfrak{w}$.
The S\={a}to weights of coordinates are  $\wgt u_{w_k} \,{=}\, -w_k$. 
The holomorphic differentials on a trigonal curve are introduced as follows,
cf.\;\eqref{GSeqTrig}:
\begin{subequations}\label{HolDiffs}
\begin{align}
&\mathfrak{C}^1:& &\rmd u_{3k-1} = \frac{x^{2m-k} \rmd x}{\partial_y f},\qquad k=1,\dots, 2\mFr,\\
&& &\rmd u_{3k-2} = \frac{y x^{m-k} \rmd x}{\partial_y f},\qquad k=1,\dots, \mFr; \notag \\
&\mathfrak{C}^2:& &\rmd u_{3k-2} = \frac{x^{2m+1-k} \rmd x}{\partial_y f},\qquad k=1,\dots, 2\mFr+1,\\
&& &\rmd u_{3k-1} = \frac{y x^{m-k} \rmd x}{\partial_y f},\qquad k=1,\dots, \mFr. \notag 
\end{align}
\end{subequations}
Let $\{\mathfrak{a}_j$, $\mathfrak{b}_j \mid j=1$, \ldots $g \}$ 
form a canonical homology basis on a curve, and $\omega$,
$\omega'$ be period matrices with entries
\begin{gather*}
\omega_{w_k,j} = \int_{\mathfrak{a}_j} \rmd u_{w_k},\qquad\qquad
\omega'_{w_k,j} = \int_{\mathfrak{b}_j} \rmd u_{w_k}.
\end{gather*}
The Jacobian variety of a curve is the factor $\Complex^g/ \mathfrak{P}$,
where $\mathfrak{P}$ is the lattice of periods $\omega$, $\omega'$.
Abel's map $\mathcal{A}:C \mapsto  \Jac(C)$ is defined by
\begin{gather}
\mathcal{A}(P) = \int_\infty^P \rmd u,
\end{gather}
where $\rmd u = (\rmd u_{w_1}, \rmd u_{w_2},\dots, \rmd u_{w_g})^t$, and
on the $n$-th symmetric power $C^{(n)}$ of a curve $C$ as 
\begin{gather}
\mathcal{A}(D) = \sum_{k=1}^n \mathcal{A}(P_k),\qquad D = \sum_{k=1}^n P_k.
\end{gather}

\subsection{Non-special divisors}\label{ss:NSpDiv}
Let $D=\sum_{k=1}^g (x_k,\,y_k)$ be a divisor, and  $l(D)$ denote the dimension
of the vector space $L(D)$ of functions having poles only at points of $D$ 
with at least that multiplicities.
Let $K$ be the canonical divisor.
A divisor $D$ is called \emph{special} if $l(K-D)>0$, and
 $l(K-D)$ is its \emph{index of speciality} $i(D)$, see \cite[p.\;296]{Hart}.
Otherwise, $D$ is \emph{non-special}. 
Evidently, such a definition is not very helpful if one needs to determine 
whether a given divisor is special or non-special.

In \cite[III.6.6., p.\,92]{Far80} the following criteria of a non-special degree $g$ divisor can be found.
The condition $i(D)=0$ gives immediately 
that the rank of the Jacobian of $\mathcal{A}(D)$ is $g$. In fact, 
the Jacobian can be replaced by a matrix $\mathfrak{R}_{g}[D]$ constructed as follows.

Let $D=\sum_{i=1}^{N} (x_i,\,y_i)$, then 
$\mathfrak{R}_{m}[D]$ is a matrix of size $N \times m$ constructed from the first $m$ 
monomials of the ordered list $\mathfrak{M}$ associated with a curve in such a way
that the $i$-th row of $\mathfrak{R}_{m}[D]$ consists of the values of these monomials at $(x_i,y_i)$.

\newtheorem*{SpecDivCrit}{Special Divisor Criteria}
\begin{SpecDivCrit}
An effective divisor $D$ on an algebraic curve is non-special if and only if 
$\rank \mathfrak{R}_{g}[D] = g$.  $D$ is special if $\rank \mathfrak{R}_{g}[D] < g$.
\end{SpecDivCrit}

In what follows we define divisors with the help of 
entire rational functions, see Section~\ref{s:RatF} for more details.
 If $\mathcal{B}$ is an entire rational function with 
a divisor $Z$ of zeros, $\deg Z=N$, then $(\mathcal{B}) = Z -N \infty$
is the principal divisor of $\mathcal{B}$.
In other words, we work with  the vector spaces $L(n \infty)$,
and $\mathcal{B}$ serves as a section in $L(N \infty)$,

The \emph{involution on a trigonal curve} connects three points with the same $x$-coordinates,
say $(a,\,b_1)$, $(a,\,b_2)$, $(a,\,b_3)$. 
These three points are zeros of the entire rational function
$\mathcal{K}(x,y) = x-a$ on the curve, and so $(\mathcal{K}) = \sum_{k=1}^3 (a,\,b_k) - 3\infty$.

\begin{lem}
Every group of three points in involution
in an effective divisor $D$ on a trigonal curve generates one linear dependent row in 
$\mathfrak{R}_g[D]$.
\end{lem}
\begin{proof}
This result is obtained by direct computations. The first $g$ monomials on a trigonal curve
are
\begin{align*}
&\mathfrak{C}^1:\quad \mathfrak{M} = \big\{1,  x, \dots, x^{\mFr-1}, x^\mFr, y, x^{\mFr+1}, 
y x, \dots, x^{2\mFr-1}, y x^{\mFr-1} \big\},  \\  
&\mathfrak{C}^2:\quad \mathfrak{M} = \big\{ 1,  x,  \dots, x^{\mFr-1},  x^\mFr,  y,  x^{\mFr+1},  
y x,  \dots, x^{2\mFr-1}, y x^{\mFr-1}, x^{2\mFr} \big\}.
\end{align*}
Note, that all monomials are linear in $y$. Suppose, $D$ contains a
group of three points connected by involution, namely: $(a,\,b_1)$, $(a,\,b_2)$, $(a,\,b_3)$.
Then  $\mathfrak{R}_g[D]$ contains the following rows
\begin{multline*}
\begin{pmatrix}
1&  a& \cdots& a^\mFr& b_1 & a^{\mFr+1}& b_1 a& \cdots& a^{2\mFr-1}& b_1 a^{\mFr-1} \\
1&  a& \cdots& a^\mFr& b_2 & a^{\mFr+1}& b_2 a& \cdots& a^{2\mFr-1}& b_2 a^{\mFr-1} \\
1&  a& \cdots& a^\mFr& b_3 & a^{\mFr+1}& b_3 a& \cdots& a^{2\mFr-1}& b_3 a^{\mFr-1} 
\end{pmatrix} 
\begin{matrix} R_1 \\ R_2 \\ R_3  \end{matrix} \\
\sim
\left( \begin{matrix}
1&  a& \cdots& a^\mFr& b_1 & a^{\mFr+1}& b_1 a& \cdots \\
0&  0& \cdots&  0& b_2-b_1 & 0 & (b_2-b_1) a& \cdots \\
0&  0& \cdots&  0& b_3-b_1 & 0 & (b_3-b_1) a& \cdots 
\end{matrix} \right. \phantom{mmmmmmmmmmm} \\
\left. \begin{matrix}
 a^{2\mFr-1}& b_1 a^{\mFr-1} \\
 0 & (b_2-b_1) a^{\mFr-1} \\
 0 & (b_3-b_1) a^{\mFr-1}
\end{matrix} \right)
\begin{matrix} R_1 \, \phantom{= R_2-R_1}  \\  R'_2  = R_2-R_1\\ R'_3 = R_3-R_1 
\end{matrix}
\end{multline*}
on a curve of type $\mathfrak{C}^1$. Evidently, rows $R'_2$ and $R'_3$ coincide 
up to a constant multiple. The same situation takes place on a curve of type $\mathfrak{C}^2$.
\end{proof}

\begin{defn}
We call a non-special divisor \emph{strictly non-special} if it contains no two points connected by involution. 
\end{defn}

\subsection{Theory of addition laws}\label{ss:AddLaws}
Here we call attention to relatively new results which put in order the theory of addition laws
on algebraic curves, see \cite{bl05}.

Let $\textsf{Y}= \Complex^n$ with coordinates $\lambda=(\lambda_k)$ be the space of parameters
of a family of algebraic curves of a fixed genus $g$, cf \eqref{CnsEq}. Let $\mathfrak{V}$
be a family of curves $C_{n,s}$ with fixed $n$ and $s$. Let $\textsf{X}$
be the universal fiber-bundle of the $g$-th symmetric power of curves from $\mathfrak{V}$
over $\textsf{Y}$:
\begin{equation*}
  0 \to \mathfrak{V}^{(g)} \to  \textsf{X} \overset{\textsf{p}}{\rightarrow}  \textsf{Y}   \to 0.
\end{equation*}
A point $\textsf{x} \in \textsf{X}$ is an unordered set of 
$g$ points $(x_i,y_i)\in \Complex^2$ of a curve $C \in \mathfrak{V}$ defined by $f(x,y,\lambda)=0$,
and $\textsf{p}(\textsf{x}) = \lambda$, where $\lambda$ are parameters of the curve $C$. 

Let $\textsf{X} \times_{\textsf{Y}} \textsf{X}$
be the direct product of  $\textsf{X}$ over $\textsf{Y}$.

\begin{prop*}[Theorem 3.2 \cite{bl05}]
The space $\textsf{X}$ together with the mapping $\textsf{p}: \textsf{X} \to \textsf{Y}$
and two structure mappings 
$\add: \textsf{X} \times_{\textsf{Y}} \textsf{X}
\to \textsf{X}$, and $\inv: \textsf{X} \to \textsf{X}$ such that
\begin{align*}
& \add(\add(\textsf{x}_1,\, \textsf{x}_2),\, \textsf{x}_3) = 
\add(\textsf{x}_1,\, \add(\textsf{x}_2,\, \textsf{x}_3)), \\
& \add(\add(\textsf{x}_1,\, \textsf{x}_2),\, \inv(\textsf{x}_2)) = \textsf{x}_1,
\end{align*}
provided $\textsf{p}(\textsf{x}_1) = \textsf{p}(\textsf{x}_2) 
= \textsf{p}(\textsf{x}_3)$,
form a commutative algebraic groupoid over $\textsf{Y}$.
\end{prop*}
Note,  $\add$ may not be defined on all
pairs $\textsf{x}_1$ and $\textsf{x}_2$ from  $\textsf{X}$.

The introduced structure mappings stand for an addition law on $\textsf{X}$,
which is realized between two degree $g$ non-special divisors.  
Mapping $\inv$ is constructed with the help of an entire rational function $\mathcal{R}_{2g}$
of weight~$2g$, and mapping $\add$ with the help of an entire rational function $\mathcal{R}_{3g}$
of weight~$3g$. For more details see \cite{bl05}.

In fact, an entire rational function $\mathcal{R}_N$ of weight $N > 2g$  
can be used to reduce a non-special divisor of degree $N-g$ to a degree $g$ one. 
The function $\mathcal{R}_N$ itself produces the inverse divisor to the reduced one.
Then a function $\mathcal{R}_{2g}$ of weight $2g$ constructed from the inverse divisor
gives the desired reduced divisor.

In this paper we propose to add divisors using such a reduction. 
And we suggest to add point by point, so only functions $\mathcal{R}_{2g}$
and $\mathcal{R}_{2g+1}$ of weights $2g$ and $2g+1$ are involved.

\section{Divisors in terms of entire rational functions}\label{s:RatF}

\subsection{Entire rational functions on a trigonal curve}
An entire rational function on a curve is a function with a finite number of zeros which are points of the curve.
The number of zeros is shown by the S\={a}to
 weight of an entire function. So a function $\mathcal{B}$ of weight $N$ has $N$ zeros, 
and  $\wgt \mathcal{B} = N$. 
The divisor of zeros $Z$ and the divisor of poles $N \infty$ of an entire function form a principal divisor: 
$(\mathcal{B}) = Z - N \infty$. 
According to the Riemann-Roch theorem, on a curve of genus $g$ 
among $N$ zeros of the entire function~$\mathcal{B}$  
only $N-g$ can be chosen arbitrarily, we denote them by~$D$. 
The remaining $g$ zeros form the divisor $D^\ast$ complement to $D$
in $Z$,  that is $D + D^\ast = Z$. We call $D^\ast$ the inverse divisor
to $D$ with respect to $\mathcal{B}$.

\begin{defn}
We call a degree $g$ non-special divisor $D^\ast$ the \emph{inverse divisor} to 
a degree $p$ non-special divisor $D$, $p\geqslant g$, if $D + D^\ast$ is
the divisor of zeros of an entire rational function of weight $g+p$.
\end{defn}

In what follows divisors $Z$, $D$, $D^\ast$
are supposed to be effective.
We use the algebra of entire rational functions to define such divisors. 
We often omit the word `rational', so 
`an entire function' stands for `an entire rational function'

An entire function on a plane algebraic curve is
a linear combination of monomials 
from the list $\mathfrak{M}$. In fact, an entire function is a polynomial in $x$ and $y$. 
We refer to  the list $\mathfrak{M}$ in order to control the weight of a function.
An entire function $\mathcal{B}$ of weight~$N$, $N\geqslant g + 1$,
contains the first $N-g+1$ elements 
from the list of monomials. Note,  there is no entire functions of the weights 
equal to elements of the Weierstrass gap sequence.

\textbf{I.} When $N > 2g+1$,  function $\mathcal{B}$ on a trigonal curve has the form 
\begin{subequations}\label{Bfunct2}
\begin{gather}\label{Bform2}
\mathcal{B}(x,y) = y^2 \beta_{2y}(x) + y\beta_y(x) + \beta_x(x),
\end{gather}
and $N \equiv \wgt \mathcal{B}$ is arbitrary within the interval.
Let $p =N-2g$, then
\begin{align} \label{Bdegs2}
&\mathfrak{C}^1:\ \deg \beta_{2y} = [(p -2)/3],& 
&\mathfrak{C}^2:\ \deg \beta_{2y} = [(p -2)/3]& \\
&\phantom{\mathfrak{C}^1:\ \,}\deg \beta_{y} = \mFr+ [(p -1)/3],& 
&\phantom{\mathfrak{C}^2:\ \,}\deg \beta_{y} = \mFr+ [p/3],& \notag\\
&\phantom{\mathfrak{C}^1:\ \,}\deg \beta_{x} = 2\mFr+ [p/3];& 
&\phantom{\mathfrak{C}^2:\ \,}\deg \beta_{x} = 2\mFr+ [(p +2)/3],& \notag
\end{align}
\end{subequations}
where $[z]$ denotes the integer part of $z$. 
\begin{theo}\label{T:BfSyst2}
Let $\mathcal{B}$ be an entire rational function of weight $N$, $N > 2g+1$, 
on a trigonal curve of genus $g$ 
defined by \eqref{C3def},
that is $\mathcal{B}$ has the form \eqref{Bfunct2} 
and its $N$ zeros are solutions of the system
\begin{gather}\label{BfEqs2}
\mathcal{B}(x,y) = 0,\qquad f(x,y)=0.
\end{gather}
Then the zeros of $\mathcal{B}$ with distinct $x$-coordinates satisfy the system
\begin{subequations}\label{Zf2gp}
\begin{gather}\label{yZf2gpSyst}
\mathcal{Z}_{N}(x) = 0,\qquad \mathcal{R}(x,y)=0,
\end{gather}
where
\begin{align} \label{yDef2gp}
& \mathcal{R}(x,y) = y \rho_y(x) + \rho_x(x),\\
\begin{split} \notag
 &\qquad \rho_y(x) = \beta_{y}^2 - \beta_{2y}\beta_x
 - \beta_{2y} \beta_{y}\mathcal{T} + \beta_{2y}^2  \mathcal{Q},\\
 &\qquad \rho_x(x) = \beta_{y}\beta_x - \beta_{2y}\beta_x \mathcal{T}  
+  \beta_{2y}^2  \mathcal{P},
\end{split}\\
&\mathcal{Z}_{N}(x) = \Resultant_y (\mathcal{B},f)  \label{Z2gp}\\
 &\phantom{\mathcal{Z}_{g+p}(x)} = \beta_x^3 - \beta_x^2 \beta_y \mathcal{T} 
 + \beta_x^2 \beta_{2y} (\mathcal{T} ^2 - 2\mathcal{Q})
 + \beta_x \beta_y^2 \mathcal{Q} + \beta_x \beta_y \beta_{2y} (3\mathcal{P} - \mathcal{Q}\mathcal{T}) \notag \\
 &\phantom{\mathcal{Z}_{g+p}(x)}  \quad  + \beta_x \beta_{2y}^2 (\mathcal{Q}^2 - 2\mathcal{P}\mathcal{T})
 - \beta_y \beta_{2y}^2 \mathcal{P}\mathcal{Q} + \beta_y^2 \beta_{2y} \mathcal{P}\mathcal{T}
 - \beta_y^3 \mathcal{P} + \beta_{2y}^3 \mathcal{P}^2. \notag
\end{align}
\end{subequations}
\end{theo}
A proof is the matter of straightforward computation. 

\medskip

\textbf{II.} When $g+\mFr < N \leqslant 2g+1$, an entire function $\mathcal{B}$ has the form 
\begin{subequations}\label{Bfunct1}
\begin{gather}\label{Bform1}
\mathcal{B}(x,y) = y \beta_y(x) + \beta_x(x),
\end{gather}
and $N \equiv \wgt \mathcal{B} \neq 3(\mFr+k) -1$, $k=1$, \ldots, $\mFr$,
on curves from $\mathfrak{C}^1$, or $N \equiv \wgt \mathcal{B} \neq 3(\mFr+k) -2$, $k=1$, \ldots, $\mFr$,
on curves from $\mathfrak{C}^2$.
Let $p=N-g$, then
\begin{align} \label{Bdegs1}
&\mathfrak{C}^1:\ \deg \beta_{y} =  [(p-1)/3],& 
&\mathfrak{C}^2:\ \deg \beta_{y} =  [(p-1)/3]& \\
&\phantom{\mathfrak{C}^1:\ \,}\deg \beta_{x} = \mFr+ [p/3];& 
&\phantom{\mathfrak{C}^2:\ \,}\deg \beta_{x} = \mFr+ [(p+1)/3].& \notag
\end{align}
\end{subequations}
\begin{theo}\label{T:IfSyst1}
Let $\mathcal{B}$ be an entire rational function of weight  $N$, $g+\mFr < N \leqslant 2g+1$, 
 on a trigonal curve of genus $g$
defined by \eqref{C3def}, that is $\mathcal{B}$ has the form \eqref{Bfunct1}
and its $N$ zeros are solutions of the system
\begin{gather}\label{IfEqs1}
\mathcal{B}(x,y) = 0,\qquad f(x,y)=0.
\end{gather} 
Then the zeros of $\mathcal{B}$ with distinct $x$-coordinates satisfy the system
\begin{subequations}\label{Zf2g1}
\begin{gather}\label{yZf2g1Syst}
 \mathcal{Z}_{N}(x) = 0, \qquad \mathcal{B}(x,y) = 0,
\end{gather}
where
\begin{align}\label{Zgp}
& \mathcal{Z}_{N}(x) = 
\Resultant_y (\mathcal{B},f)  = \beta_y^3 f\big(x,- \beta_x/ \beta_y\big) \\
&\phantom{\mathcal{Z}_{g+p}(x)} 
 = -\beta_x^3 + \beta_x^2 \beta_y \mathcal{T}
 - \beta_x \beta_y^2 \mathcal{Q} + \beta_y^3 \mathcal{P}. \notag
\end{align}
\end{subequations}
\end{theo}
A proof is the matter of straightforward computation. 

\medskip
\textbf{III.} 
\begin{theo}\label{T:IfSyst0}
Let $\mathcal{B}$ be an entire rational function of weight  $N$, $N \leqslant g+\mFr $, 
 on a trigonal curve of genus $g$ defined by \eqref{C3def}. Then 
\begin{gather}\label{Bform0}
\mathcal{B}(x,y) = \beta_x(x),
\end{gather}
and $N \equiv \wgt \mathcal{B} = 3 \deg \beta_x$. 
The divisor of zeros of  $\mathcal{B}$ consists of
$N/3$ groups of three points in involution, the $x$-coordinate in each group 
is a zero of $\beta_x$.
\end{theo}

In the case of a trigonal curve, the list  $\mathfrak{M}$ is given by \eqref{Monoms}. With 
a degree $g$ divisor $D$ we construct 
\begin{gather} \label{RgMatr}
\begin{split}
\mathfrak{C}^1:\quad
&\mathfrak{R}_{g}[D] =  
\big\| \!\!\!
\begin{array}{rccccl} 
 \{1 & x_k \ \dots \ x_k^{\mFr} & y_k & x_k^{\mFr+1} & y_k x_k \ \dots \  
 x_k^{2\mFr-1} & y_k x_k^{\mFr-1} \}_{k=1}^g
 \end{array} \!\!\! \big\|, \\
\mathfrak{C}^2:\quad
&\mathfrak{R}_{g} [D] =  
\big\| \!\!\!
\begin{array}{rccccl} 
 \{1 & x_k \ \dots \ x_k^{\mFr} & y_k & x_k^{\mFr+1} & y_k x_k \ \dots \  
 y_k x_k^{\mFr-1}  & x_k^{2\mFr} \}_{k=1}^g
 \end{array} \!\!\! \big\|.   
 \end{split}
 \end{gather}

\begin{prop}\label{P:ERFdetF}
An entire rational function $\mathcal{B}$ of weight $g+p$, $p \geqslant \mFr $, 
on a trigonal curve \eqref{C3def} of genus~$g$
is constructed from a collection of $p$ points with the help of the determinant formula 
based on the first $p+1$ monomials from the list \eqref{Monoms}. 

Each selection of $p$ points from the divisor of zeros of $\mathcal{B}$ produces  
the same function $\mathcal{B}$ up to a constant multiple.
\end{prop}
\begin{proof}
Let $D_p = \sum_{k=1}^{p} (x_k,y_k)$ be a collection of distinct points of \eqref{C3def} 
without points connected by involution.
The determinant formula is obtained as follows.
The first $p+1$ monomials from the list \eqref{Monoms} form the first row
of a square matrix $\overline{\mathfrak{R}}[D_p]$ of size $p+1$.
Values of these monomials at $p$ points of $D_p$ form the next $p$ rows,
which is  $\mathfrak{R}_{p+1}[D_p] $.
Then $\mathcal{B}(x,y) = \det \overline{\mathfrak{R}}[D_p]$. 
The determinant formulas for entire rational functions of weights $2g$ and $2g+1$
on trigonal curves 
are described in detail in subsections~\ref{ss:RatF2g} and \ref{ss:RatF2g1}.

Suppose, $\mathcal{B}$ and $\widetilde{\mathcal{B}}$ are entire rational functions of the same weight,
with the same divisor of zeros $Z$. 
Let $\mathcal{B}$ be constructed from a selection of $p$ points $D_p \subset Z$, 
and $\widetilde{\mathcal{B}}$  from another selection of $p$ points 
$\widetilde{D}_p \subset Z$. If $p>g+1$ the both functions have the form \eqref{Bfunct2}
and polynomials $\beta_{2y}$, ${\beta}_{y}$, ${\beta}_{x}$ from $\mathcal{B}$
have the same degrees, respectively, as polynomials $\widetilde{\beta}_{2y}$, $\widetilde{\beta}_{y}$, 
$\widetilde{\beta}_{x}$ from $\widetilde{\mathcal{B}}$.
By equating resultants \eqref{Z2gp} 
\begin{equation*}
\Resultant_y (\mathcal{B},f)  = c \Resultant_y (\widetilde{\mathcal{B}},f) 
\end{equation*}
with an arbitrary constant multiple $c$,
and taking into account that parameters of a curve $f(x,y)=0$ may vary,
we find the following relations:
\begin{align*}
&\beta_{2y}^3 = c \widetilde{\beta}_{2y}^3,&  
&\beta_{2y} \beta_{y} \beta_{x} = c \widetilde{\beta}_{2y} \widetilde{\beta}_{y} \widetilde{\beta}_{x}&
&\beta_{y}^3 - 3 \beta_{2y} \beta_{y} \beta_{x} 
= c \big(\widetilde{\beta}_{y}^3 - 3 \widetilde{\beta}_{2y} \widetilde{\beta}_{y} \widetilde{\beta}_{x}\big),&\\
&\beta_{2y}^2 \beta_{y} = c \widetilde{\beta}_{2y}^2 \widetilde{\beta}_{y},& 
&\beta_{2y}^2 \beta_{x} =c  \widetilde{\beta}_{2y}^2 \widetilde{\beta}_{x},& 
&\beta_{2y}  \beta_{y}^2 - 2 \beta_{2y}^2 \beta_{x} 
= c \big(\widetilde{\beta}_{2y} \widetilde{\beta}_{y}^2  - 2\widetilde{\beta}_{2y}^2 \widetilde{\beta}_{x} \big),& \\
&\beta_{y} \beta_{x}^2 = c \widetilde{\beta}_{y} \widetilde{\beta}_{x}^2,& 
&\beta_{2y} \beta_{x}^2 = c \widetilde{\beta}_{2y} \widetilde{\beta}_{x}^2,&
&\beta_{y}^2 \beta_{x} - 2 \beta_{2y} \beta_{x}^2 
= c \big(\widetilde{\beta}_{y}^2 \widetilde{\beta}_{x} - 2\widetilde{\beta}_{2y} \widetilde{\beta}_{x}^2 \big),& \\
&\beta_{x}^3 = c \widetilde{\beta}_{x}^3.&
\end{align*}
The equalities are true for an arbitrary choice of the argument $x$ of the polynomials. Thus, the only solution is
\begin{equation*}
\beta_{2y} = \epsilon \widetilde{\beta}_{2y},\qquad
\beta_{y} = \epsilon \widetilde{\beta}_{y},\qquad
\beta_{x} = \epsilon \widetilde{\beta}_{x},\qquad \epsilon^3 = c.
\end{equation*}
Assigning $\beta_{2y}(x) \equiv 0$, we obtain the same result in the case of $\mFr < p \leqslant g + 1$,
when the functions $\mathcal{B}$ and $\widetilde{\mathcal{B}}$ have the form \eqref{Bfunct1}.
Therefore, the functions $\mathcal{B}$ and $\widetilde{\mathcal{B}}$ coincide up to a constant multiple. 
\end{proof}

\begin{rem}
Note, that an entire rational function $\mathcal{B}$ of weight $N$ such that $g + \mFr < N \leqslant 2g+1$ 
is linear in~$y$. Thus, divisors with two points in involution do not arise among solutions of the system \eqref{IfEqs1}.
Moreover, if such a function $\mathcal{B}$ is constructed from a divisor $D$ of degree $N-g$
which contains two points in involution, then the third point from this involution group belongs to the divisor $D^\ast$
inverse to $D$ with respect to~$\mathcal{B}$. This is seen from the direct computations made for 
entire functions of weights $2g$ and $2g+1$ below.
\end{rem}

Entire functions of weights $2g$ and $2g+1$ play the key role in the reduction algorithm.
We describe them in more detail.

\subsection{Entire function of order $2g$}\label{ss:RatF2g}
Let  $D_g = \sum_{k=1}^g (x_k,y_k)$ be a collection of distinct points 
on a trigonal curve of genus~$g$, no points connected by involution.
According to Proposition~\ref{P:ERFdetF}, 
an entire function $\mathcal{I}$ of weight $2g$ is constructed from $D_g$
with the help of the determinant formula as follows. 
A square matrix $\overline{\mathfrak{R}}[D_g]$
of size $g+1$ has the form
\begin{subequations}\label{R2gMatr}
\begin{gather}\label{R2gMatr1}
\begin{split}
\mathfrak{C}^1:\quad
&\overline{\mathfrak{R}} [D_g] =  
\bigg\| \!\!\!
\begin{array}{rccccl} 1 & x \ \ \dots \ \ x^{\mFr} & y & x^{\mFr+1} & y x \ \  \dots \ \  y x^{\mFr-1} & x^{2\mFr}  \\ 
 \{1 & x_k \ \dots \ x_k^{\mFr} & y_k & x_k^{\mFr+1} & y_k x_k \ \dots \  y_k x_k^{\mFr-1} & x_k^{2\mFr}\}_{k=1}^g
 \end{array} \!\!\! \bigg\|,\\
\mathfrak{C}^2:\quad
&\overline{\mathfrak{R}} [D_g]=  
\bigg\| \!\!\!
\begin{array}{rccccl} 1 & x \ \ \dots \ \ x^{\mFr} & y & x^{\mFr+1} & y x \ \  \dots \ \  x^{2\mFr}  & y x^\mFr \\ 
 \{1 & x_k \ \dots \ x_k^{\mFr} & y_k & x_k^{\mFr+1} & y_k x_k \ \dots \   x_k^{2\mFr}& y_k x_k^{\mFr} \}_{k=1}^g
 \end{array} \!\!\! \bigg\|.
 \end{split}
\end{gather}
The entire function $\mathcal{I}$ of weight $2g$ vanishing on $D_g$ is
\begin{equation}\label{IDetR}
 \mathcal{I}(x,y) = \det \overline{\mathfrak{R}}[D_g]. 
\end{equation}
\end{subequations}

\begin{rem}\label{R:repPs}
Let  $D_g$ contain a pair of repeated points, say $(x_g,y_g)=(x_1,y_1)$.
Then the last row of $\overline{\mathfrak{R}} [D_g]$, 
corresponding to the point $(x_g,y_g)$, should be replaced by the limit as $(x_g,y_g)\to (x_1,y_1)$ 
of the total derivative with respect to $x_g$.  That is in \eqref{R2gMatr1} the last row has the form
\begin{multline*}
\mathfrak{C}^1:\quad  \{0,\  1, \ \dots, \ \mFr x_1^{\mFr-1}, \ y'_1, \ (\mFr+1) x_1^{\mFr},\  y_1 +  x_1 y'_1, \ \dots, \\   
 (\mFr-1) x_1^{\mFr-2} y_1  + x_1^{\mFr-1}  y'_1,\ 2\mFr x_1^{2\mFr-1} \},
\end{multline*}
where $y'_1 =  \rmd y_1 / \rmd x_1 = - \partial_{x_1} f(x_1,y_1)/\partial_{y_1} f(x_1,y_1)$.
If the multiplicity of a point $(x_k,\, y_k)$ is greater than $2$, say $n$, then $n-1$ rows,
corresponding to the repeated point $(x_k,\, y_k)$, are replaced by the total derivatives of orders
up to $n-1$.
\end{rem}

We call $\mathcal{I}$ the \emph{minimal function} defining a degree $g$ divisor $D_g$.
It has the form
\begin{gather}\label{IDef3}
\mathcal{I}(x,y) = y \alpha_{y} (x) + \alpha_{x} (x) ,\\ \notag
\begin{aligned}
&\mathfrak{C}^1:& &\deg (\alpha_{y}) = \mFr-1,& &\deg (\alpha_{x}) = 2\mFr,& \\  
&\mathfrak{C}^2:& &\deg (\alpha_{y}) = \mFr,& &\deg (\alpha_{x}) = 2\mFr.&
\end{aligned}
\end{gather}

In what follows we need the entire functions $\mathcal{I}$ on $C_{3,4}$ 
constructed from (i)~three distinct points, (ii)~three equal points.
\begin{exam}\label{E:I3AbPC1}
Let $D_3 = \sum_{k=1}^3 (a_k,\,b_k)$ on $C_{3,4}$. Then
\begin{subequations}\label{IDetRC34}
\begin{gather}
\overline{\mathfrak{R}} [D_3] = 
\begin{pmatrix}
1 & x & y & x^2 \\
1 & a_1 & b_1 & a_1^2 \\
1 & a_2 & b_2 & a_2^2 \\
1 & a_3 & b_3 & a_3^2 \\
\end{pmatrix},
\end{gather}
and
\begin{multline}
\mathcal{I}(x,y) = \det \overline{\mathfrak{R}} [D_3]
= y (a_2-a_1)(a_3-a_1)(a_3-a_2) \\
- x^2 \big( b_1 (a_3-a_2) + b_2 (a_1 - a_3) +  b_3 (a_2-a_1) \big) \\
+ x \big(b_1 (a_3^2 - a_2^2) + b_2 (a_1^2 - a_3^2) +  b_3 (a_2^2 - a_1^2)\big) \\
- \big(b_1 a_2 a_3 (a_3-a_2) + b_2 a_1 a_3 (a_1-a_3) + b_3 a_1 a_2 (a_2-a_1)\big).
\end{multline}
\end{subequations}
\end{exam}

\begin{exam}\label{E:I3PC1}
On the curve \eqref{C34Ex1} let $D_3=(-1,\,5)+(3,\,1)+(4,-3)$. 
By \eqref{IDetRC34} we have
\begin{equation*}
\mathcal{I}(x,y)  = 4 (3x^2 + 5y - x - 29).
\end{equation*}
In the further computations we always cancel the common numerical factor.

Note, that the divisor $D_3$ is unambiguously defined by the system of the polynomial 
$\mathcal{H}(x)= (x+1)(x-3)(x-4) $ 
and the obtained minimal function $\mathcal{I}$, namely:
\begin{subequations}
\begin{align}\label{C1Div3}
&x^3 - 6x^2 + 5x + 12 = 0, \\
&3x^2 + 5y - x - 29 = 0.  \label{C1Div3}
\end{align}
\end{subequations}
\end{exam}

\begin{exam}\label{E:I3SamePAbC1}
On a trigonal curve of the form \eqref{C3def} we find
\begin{align*}
&y' = - \frac{\partial_x f}{\partial_y f} 
= - \frac{y^2  \mathcal{T}'(x)  + y \mathcal{Q}'(x) + \mathcal{P}'(x)}
{3 y^2 +2 y  \mathcal{T}(x) +  \mathcal{Q}(x)}& \\
&y'' = - \frac{\partial^2_x f}{\partial_y f} + 2  \frac{(\partial_x f) \partial_y \partial_x f}{(\partial_y f)^2}
- \frac{(\partial_x f)^2 \partial^2_y f}{(\partial_y f)^3}& \\
&\quad = - \frac{y^2  \mathcal{T}''(x)  + y \mathcal{Q}''(x) + \mathcal{P}''(x)}
{3 y^2 +2 y  \mathcal{T}(x) +  \mathcal{Q}(x)}\\
&\quad\quad + 2  \frac{\big(y^2  \mathcal{T}'(x)  + y \mathcal{Q}'(x) + \mathcal{P}'(x) \big)
\big(2 y  \mathcal{T}'(x) +  \mathcal{Q}'(x) \big)}
{\big( 3 y^2 + 2 y  \mathcal{T}(x) +  \mathcal{Q}(x) \big)^2} \\
&\quad\quad - \frac{\big(y^2  \mathcal{T}'(x)  + y \mathcal{Q}'(x) + \mathcal{P}'(x) \big)^2
\big(6y + 2 \mathcal{T}(x) \big)}
{\big(3 y^2 +2 y  \mathcal{T}(x) +  \mathcal{Q}(x) \big)^3}.
\end{align*}
In the case of $C_{3,4}$ defined by \eqref{C34eq1}, we obtain
\begin{align*}
&y'  = - \frac{y^2 \lambda_1 + y \big( 2 \lambda_8 x + \lambda_5 \big)  
+ 4 x^3 +3  \lambda_3 x^2 + 2 \lambda_6 x + \lambda_9}
{3y^2 + 2y \big( \lambda_1 x + \lambda_4 \big) + \lambda_8 x^2 + \lambda_5 x + \lambda_8 } \\
&y'' = - \frac{  2 \lambda_8 y + 12 x^2 + 6 \lambda_3 x + 2 \lambda_6 }
{3y^2 + 2y \big( \lambda_1 x + \lambda_4 \big) + \lambda_8 x^2 + \lambda_5 x + \lambda_8 } \\
&\quad \quad + 2  \frac{\big(y^2 \lambda_1 + y \big( 2 \lambda_8 x + \lambda_5 \big)  
+ 4 x^3 +3  \lambda_3 x^2 + 2 \lambda_6 x + \lambda_9\big)
\big(2\lambda_1 y + 2 \lambda_8 x + \lambda_5 \big)}
{\big(3y^2 + 2y \big( \lambda_1 x + \lambda_4 \big) + \lambda_8 x^2 + \lambda_5 x + \lambda_8 \big)^2} \\
&\quad\quad - \frac{\big(y^2 \lambda_1 + y \big( 2 \lambda_8 x + \lambda_5 \big)  
+ 4 x^3 +3  \lambda_3 x^2 + 2 \lambda_6 x + \lambda_9\big)^2
\big(6y + 2\lambda_1 x + 2 \lambda_4\big)}
{\big(3y^2 + 2y \big( \lambda_1 x + \lambda_4 \big) + \lambda_8 x^2 + \lambda_5 x + \lambda_8 \big)^3}.
\end{align*}

Let $D_3 = 3 (a,\,b)$. Then
\begin{subequations}\label{IDetRC34SameP}
\begin{gather}
\overline{\mathfrak{R}} [D_3] = 
\begin{pmatrix}
1 & x & y & x^2 \\
1 & a & b & a^2 \\
0 & 1 & y'(a,b) & 2a \\
0 & 0 & y''(a,b) & 2
\end{pmatrix},
\end{gather}
and
\begin{gather}
\mathcal{I}(x,y) = \det \overline{\mathfrak{R}} [D_3]
= 2 (y-b) - y''(a,b) (x-a)^2 - 2 y'(a,b) (x-a).
\end{gather}
\end{subequations}
\end{exam}

\begin{exam}\label{E:I3SamePC1}
On the curve \eqref{C34Ex1}  let 
$D_3=3(2,\,-1)$.  By \eqref{IDetRC34SameP} we have
\begin{gather}
\mathcal{I}(x,y) = 2\,162\, x^2 + 243\, y - 10\,457\, x + 12\,509.
\end{gather}
The given divisor $D_3$ is unambiguously defined by the system of the polynomial 
$\mathcal{H}(x)= (x-2)^3$ 
and the obtained minimal function $\mathcal{I}$.

\end{exam}

\begin{lem}
Let $D_g = \sum_{k=1}^g (x_k,\,y_k)$ be a strictly non-special divisor on a trigonal curve
of genus $g$,
all points are distinct, $W_g$ denote the determinant of the Vandermonde 
matrix constructed from $\{x_k\}_{k=1}^g$,
$\mathcal{I}$ of the form \eqref{IDef3} be the minimal function defining $D_g$,
and $\dFr_{\mathcal{I}} =\deg (\alpha_{y})$.
Then
\begin{subequations}\label{alphaCoefs}
\begin{align}
&\alpha_y(x) = \frac{W_g}{\dFr_{\mathcal{I}} !} 
\sum_{\substack{k_1,\dots,k_{\dFr_{\mathcal{I}}} = 1 \\ k_1 \neq \dots \neq k_{\dFr_{\mathcal{I}}}}}^{g} 
\bigg( \prod_{l=1}^{\dFr_{\mathcal{I}}} y_{k_l}\bigg) L^{(g)}_{k_1,\dots k_{\dFr_{\mathcal{I}}}}(x),  \\
&\alpha_x(x) = \frac{-W_g}{(\dFr_{\mathcal{I}} +1)!} 
\sum_{\substack{k_1,\dots,k_{\dFr_{\mathcal{I}} +1} = 1 \\ k_1 \neq \dots \neq k_{\dFr_{\mathcal{I}} +1}}}^{g} 
\bigg(\prod_{l=1}^{\dFr_{\mathcal{I}} +1}  y_{k_l}\bigg) M^{(g)}_{k_1,\dots k_{\dFr_{\mathcal{I}} +1}}(x),
\end{align}
\end{subequations}
where 
\begin{subequations}\label{LMdefs}
\begin{align}
&L^{(g)}_{k_1,\dots k_\dFr} (x) = \prod_{l=1}^{\dFr}
 \frac{ (x-x_{k_l})}{\prod_{\substack{i=1\\ i\neq  k_1,\dots k_\dFr}}^{g} (x_i - x_{k_l})}, \\
&M^{(g)}_{k_1,\dots,k_{\dFr}}(x) = \prod_{\substack{i=1\\ i\neq k_1,\dots k_{\dFr}}}^{g}
 \frac{(x_i - x)}{\prod_{l=1}^{\dFr} (x_i - x_{k_l})}.
\end{align}
\end{subequations}
\end{lem}
\begin{proof}
The representation \eqref{alphaCoefs} is obtained from the 
model of $\mathcal{I}$ given in \eqref{R2gMatr}--\eqref{IDef3} 
by straightforward computations.
\end{proof}

\begin{lem}\label{L:IZfactorize}
Suppose, a degree $g$ divisor $D_g = \sum_{k=1}^g (x_k,\,y_k)$ on a trigonal curve
of genus $g$
contains two points connected by involution, say $x_1=x_2=a$, $y_1=b_1$, $y_2=b_2$,
and $\mathcal{I}$ is constructed from $D_g$ by the determinant formula \eqref{R2gMatr}.
Then 
\begin{gather}\label{IZfactorize}
\begin{split}
\mathcal{I}(x,y) &= (b_1 - b_2) (x-a) \big(\hat{\alpha}_y(x) y + \hat{\alpha}_x(x) \big),\\
\mathcal{Z}_{2g}(x) &= (b_1 - b_2)^3 (x-a)^3 \Big( \hat{\alpha}_x(x)^3 
+ \hat{\alpha}_x(x)^2 \hat{\alpha}_y(x) \mathcal{T}(x) \\
&\quad - \hat{\alpha}_x(x) \hat{\alpha}_y(x)^2 \mathcal{Q}(x) + \hat{\alpha}_y(x)^3 \mathcal{P}(x) \Big),
\end{split}
\intertext{and}
\deg \hat{\alpha}_y = \deg \alpha_y -1,\qquad  \deg \hat{\alpha}_x = \deg \alpha_x -1. \notag
\end{gather}
This means that the third point $(a,b_3)$ connected by involution to the two mentioned above
  is located in the inverse divisor $D_g^\ast = Z - D_g \geqslant 0$,
where $Z$ denotes the divisor of zeros of $\mathcal{I}$.
\end{lem}
\begin{proof}
Consider a trigonal curve of the type $\mathfrak{C}^1$.
According to \eqref{R2gMatr}, we have
\begin{multline*}
\mathcal{I}(x,y) = 
\begin{vmatrix}
1&  x& \cdots& x^\mFr& y& x^{\mFr+1}& y x & \cdots& y x^{\mFr-1} & x^{2\mFr}  \\
1&  a& \cdots& a^\mFr& b_1 & a^{\mFr+1}& b_1 a& \cdots& b_1 a^{\mFr-1} & a^{2\mFr} \\
1&  a& \cdots& a^\mFr& b_2 & a^{\mFr+1}& b_2 a& \cdots& b_2 a^{\mFr-1} & a^{2\mFr} \\
\vdots & \vdots &  & \vdots & \vdots & \vdots & \vdots &  & \vdots & \vdots \\
1&  x_i & \cdots& x_i^\mFr& y_i& x_i^{\mFr+1}& y_i x_i & \cdots& y_i x_i^{\mFr-1} & x_i^{2\mFr} \\
\vdots & \vdots &  & \vdots & \vdots & \vdots & \vdots & \vdots & \vdots & \vdots  
\end{vmatrix} 
\begin{matrix} R_0 \\ R_1 \\ R_2 \\ \vdots \\ R_i \\ \vdots \end{matrix} \\
= \left| \begin{matrix}
0&  x-a& \cdots& x^\mFr - a^\mFr & y - b_1 & x^{\mFr+1} - a^{\mFr+1}  \\
1&  a& \cdots& a^\mFr& b_1 & a^{\mFr+1} \\
0&  0& \cdots&  0& b_2 - b_1 & 0  \\
\vdots & \vdots &  & \vdots & \vdots & \vdots \\
0&  x_i-a& \dots& x_i^\mFr - a^\mFr & y_i - b_1 & x_i^{\mFr+1} - a^{\mFr+1} \\
\vdots & \vdots &  & \vdots & \vdots & \vdots 
\end{matrix} \right. \\
\left. \begin{matrix}
 y x - b_1 a & \dots& y x^{\mFr-1} -  b_1 a^{\mFr-1} & x^{2\mFr} - a^{2\mFr} \\
 b_1 a& \dots& b_1 a^{\mFr-1} & a^{2\mFr} \\
 (b_2 - b_1) a& \dots& (b_2 - b_1) a^{\mFr-1} & 0 \\
\vdots &  & \vdots & \vdots  \\
y_i x_i - b_1 a & \dots& y_i x_i^{\mFr-1} -  b_1 a^{\mFr-1} & x_i^{2\mFr} - a^{2\mFr} \\
\vdots &  & \vdots & \vdots  
\end{matrix} \right|
\begin{matrix} R'_0 = R_0-R_1 \\ R_1 \phantom{= R_2 - R_1\ }  \\ R'_2 = R_2-R_1 \\ \vdots \\ 
R'_i = R_i - R_1 \\ \vdots  \end{matrix}.
\end{multline*}
Evidently, the substitution $x\to a$ leads to coinciding the rows $R'_0$ and $R'_2$ 
up to a constant multiple. That is $\mathcal{I}(a,y) = 0$, and so $\alpha_x$ and $\alpha_y$  
contain $(x-a)$ as a factor.
The factor $(b_1 - b_2)$ comes from $R_2'$.
Therefore, $\mathcal{I}$ and $\mathcal{Z}_{2g}$ factorize as 
in \eqref{IZfactorize}.

The case of a trigonal curve of the type $\mathfrak{C}^2$ is similar.
\end{proof}

\subsection{Entire function of order $2g+1$}\label{ss:RatF2g1}
Let  $D_{g+1} = \sum_{k=1}^{g+1} (x_k,y_k)$ be a collection of distinct points 
on a trigonal curve of genus~$g$, no points connected by involution. 
According to Proposition~\ref{P:ERFdetF}, 
an entire function $\mathcal{G}$ of weight $2g+1$ is constructed from $D_{g+1}$
with the help of the determinant formula as follows. 
A square matrix $\overline{\mathfrak{R}}_{g+2}[D_{g+1}]$ of size $g+2$ has the form
 \begin{subequations}\label{R2g1Matr}
\begin{gather} \label{R2g1Matr1}
\begin{split}
\mathfrak{C}^1:\quad 
&\overline{\mathfrak{R}}_{g+2} [D_{g+1}] =  \bigg\|\!\!\!
\begin{array}{rccccl} 1 & x \ \  \dots \ \  x^{\mFr} & y & x^{\mFr+1} & y x \ \  \dots \ \   x^{2\mFr} & y x^{\mFr} \\ 
 \{1 & x_k \ \dots \ x_k^{\mFr} & y_k & x_k^{\mFr+1} & y_k x_k \ \dots \  x_k^{2\mFr} & y_k x_k^{\mFr}\}_{k=1}^{g+1}
 \ \end{array}\!\!\!\!\! \bigg\|, \\
\mathfrak{C}^2:\quad 
&\overline{\mathfrak{R}}_{g+2} [D_{g+1}] =  \bigg\|\!\!\!
\begin{array}{rccccl} 1 & x \ \  \dots \ \  x^{\mFr} & y & x^{\mFr+1} & y x \ \  \dots \ \   y x^{\mFr} &  x^{2\mFr+1} \\ 
 \{1 & x_k \ \dots \ x_k^{\mFr} & y_k & x_k^{\mFr+1} & y_k x_k \ \dots \  y_k x_k^{\mFr} & x_k^{2\mFr+1} \}_{k=1}^{g+1}
 \ \end{array}\!\!\!\!\! \bigg\|.
 \end{split}
 \end{gather}
The entire function $\mathcal{G}$ of weight $2g+1$ vanishing on $D_{g+1}$ is
\begin{gather}\label{GdetR2g1}
 \mathcal{G}(x,y) = \det \overline{\mathfrak{R}}_{g+2}[D_{g+1}].
\end{gather}
\end{subequations} 

In the case of repeated points in $D_{g+1}$,  see Remark~\ref{R:repPs}.

We call $\mathcal{G}$ the \emph{minimal function} defining a degree $g+1$ divisor $D_{g+1}$. 
It has the form
\begin{gather}\label{GDef3}
\mathcal{G}(x,y) = y \gamma_{y} (x) + \gamma_{x} (x) ,\\ \notag
\begin{aligned}
&\mathfrak{C}^1:& &\deg \gamma_{y} = \mFr,& &\deg  \gamma_{x} = 2\mFr,& \\  
&\mathfrak{C}^2:& &\deg \gamma_{y} = \mFr,& &\deg  \gamma_{x} = 2\mFr+1.&
\end{aligned}
\end{gather}

In what follows we need the entire functions $\mathcal{G}$ on $C_{3,4}$ 
constructed from (i)~four points two of which are equal, (ii)~four equal points.

\begin{exam}\label{E:I2Same2DistPAbC1}
Let $D_4= 2(a_1,\,b_1) + (a_2,\, b_2) + (a_3,\, b_3)$  on $C_{3,4}$. Then
\begin{subequations}\label{IDetRC34Same2Dist2P}
\begin{gather}
\overline{\mathfrak{R}} [D_4] = 
\begin{pmatrix}
1 & x & y & x^2 & y x \\
1 & a_1 & b_1 & a_1^2 & b_1 a_1 \\
0 & 1 & y'(a_1,b_1) & 2a_1 &  a_1 y'(a_1,b_1) + b_1\\
1 & a_2 & b_2 & a_2^2 & b_2 a_2  \\
1 & a_3 & b_3 & a_3^2 & b_3 a_3 
\end{pmatrix},
\end{gather}
and 
\begin{multline}
\mathcal{G}(x,y) = \det \overline{\mathfrak{R}} [D_4]
=  y'(a_1,b_1) (a_1-a_2)(a_1-a_3) (x-a_1)  \Big( a_2  (y-b_3) \\ 
- a_3 (y-b_2) - x (b_2-b_3) \Big) + (a_2-a_3) (b_1-b_2)(b_1-b_3) (x-a_1)^2 \\
+ (y - b_1) (x-a_1) \Big( (b_2-b_3)(a_1-a_2)(a_1-a_3) - (b_1-b_2)(a_2-a_3)(a_1-a_3) \\
- (b_1-b_3)(a_2-a_3)(a_1-a_2) \Big) 
- (a_1-a_2)(a_1-a_3) (y - b_1) \Big( b_1 (a_2 - a_3) \\
+ b_2 (a_3 - a_1) + b_3 (a_1 - a_2) \Big).
\end{multline}
\end{subequations}
\end{exam}

\begin{exam}\label{E:I2Same2DistPC1}
On the curve \eqref{C34Ex1} let  $D_4=2(-1,\,5)+(3,\,1)+(4,\,-3)$.  
By \eqref{IDetRC34Same2Dist2P} we have
\begin{gather}\label{GSame2Dist2P}
\mathcal{G}(x,y) = 193 x y + 898 x^2 + 403 y - 1\,779 x - 3\,727.
\end{gather}
The given divisor $D_4$ is unambiguously defined by the system of the polynomial 
$\mathcal{F}(x)= (x+1)^2(x-3)(x-4)$ 
and the obtained minimal function $\mathcal{G}$.
\end{exam}

\begin{exam}\label{E:I4SamePAbC1}
Let $D_4 = 4 (a,\,b)$  on $C_{3,4}$. Then
\begin{subequations}\label{IDetRC34Same4P}
\begin{gather}
\overline{\mathfrak{R}} [D_4] = 
\begin{pmatrix}
1 & x & y & x^2 & y x \\
1 & a & b & a^2 & b a \\
0 & 1 & y'(a,b) & 2a &  a y'(a,b) + b\\
0 & 0 & y''(a,b) & 2 & a y''(a,b) + 2 y'(a,b) \\
0 & 0 & y'''(a,b) & 0 & a y'''(a,b) + 3 y''(a,b) 
\end{pmatrix},
\end{gather}
and
\begin{multline}
\mathcal{G}(x,y) = \det \overline{\mathfrak{R}} [D_4]
= - 2 y'''(a,b) \Big( (y-b) - y'(a,b) (x-a) \Big) (x-a) \\ 
- 3 y''(a,b)^2 (x-a)^2 + 6 y''(a,b) \Big( (y-b) - y'(a,b) (x-a) \Big).
\end{multline}
\end{subequations}
\end{exam}

\begin{exam}\label{E:I4SamePC1}
On the curve \eqref{C34Ex1} let
$D_4 = 4(2,\,-1)$.  By \eqref{IDetRC34Same4P} we have
\begin{gather}\label{GSame4P}
\mathcal{G}(x,y) = - 2\,156\,157\, y x  + 9\,040\,025\, x^2 + 3\,524\,265\, y  - 32\,449\,670\, x + 27\,951\,191.
\end{gather}
The given divisor $D_4$ is unambiguously defined by the system of the polynomial 
$\mathcal{F}(x)= (x-2)^4$ 
and the obtained minimal function $\mathcal{G}$.
\end{exam}

\begin{lem}
Let $D_{g+1} = \sum_{k=1}^{g+1} (x_k,\,y_k)$ be a strictly non-special divisor on a trigonal curve
of genus $g$, all points are distinct,
$W_{g+1}$ denote the determinant of the Vandermonde matrix constructed from $\{x_k\}_{k=1}^{g+1}$, 
$\mathcal{G}$ of the form \eqref{GDef3} be the minimal function defining $D_{g+1}$,
and $\dFr_{\mathcal{G}} =\deg \gamma_{y}$. Then
\begin{subequations}\label{gammaCoefs}
\begin{align}
& \gamma_y(x) = \frac{W_{g+1}}{\dFr_{\mathcal{G}} !} 
\sum_{\substack{k_1,\dots,k_{\dFr_{\mathcal{G}}} = 1 \\ k_1 \neq \dots \neq k_{\dFr_{\mathcal{G}}}}}^{g+1} 
\bigg( \prod_{l=1}^{\dFr_{\mathcal{G}}} y_{k_l}\bigg) L^{(g+1)}_{k_1,\dots k_{\dFr_{\mathcal{G}}}}(x),  \\
& \gamma_x(x) = \frac{-W_{g+1}}{(\dFr_{\mathcal{G}} +1)!} 
\sum_{\substack{k_1,\dots,k_{\dFr_{\mathcal{G}} +1} = 1 \\ k_1 \neq \dots \neq k_{\dFr_{\mathcal{G}} +1}}}^{g+1} 
\bigg(\prod_{l=1}^{\dFr_{\mathcal{G}} +1}  y_{k_l}\bigg) M^{(g+1)}_{k_1,\dots k_{\dFr_{\mathcal{G}} +1}}(x),
\end{align}
\end{subequations}
where $L^{(g)}_{k_1,\dots k_\dFr}$ and $M^{(g)}_{k_1,\dots k_\dFr}$ are defined by \eqref{LMdefs}.
\end{lem}
\begin{proof}
The representation \eqref{gammaCoefs} is obtained from the 
model of $\mathcal{G}$ given in \eqref{R2g1Matr}--\eqref{GDef3} 
by straightforward computations.
\end{proof}

\begin{lem}\label{L:GZfactorize}
Suppose, a degree $g+1$ divisor $D_{g+1}=\sum_{k=1}^{g+1} (x_k,\,y_k)$ on a trigonal curve
of genus $g$ contains two points connected by involution, say $x_1=x_2=a$, $y_1=b_1$, $y_2=b_2$,
and $\mathcal{G}$ is constructed from $D_{g+1}$ by the determinant formula \eqref{R2g1Matr}.
Then
\begin{gather*}
\begin{split}
\mathcal{G}(x,y) &= (b_1 - b_2)(x-a) (\hat{\gamma}_y(x) y + \hat{\gamma}_x(x) \big),\\
\mathcal{Z}_{2g+1}(x) &=  (b_1 - b_2)^3 (x-a)^3 \Big( \hat{\gamma}_x(x)^3 
+ \hat{\gamma}_x(x)^2 \hat{\gamma}_y(x) \mathcal{T}(x) \\
&\quad - \hat{\gamma}_x(x) \hat{\gamma}_y(x)^2 \mathcal{Q}(x) + \hat{\gamma}_y(x)^3 \mathcal{P}(x) \Big),
\end{split}
\intertext{and}
\deg \hat{\gamma}_y = \deg \gamma_y -1,\qquad  \deg \hat{\gamma}_x = \deg \gamma_x -1.
\end{gather*}
This means that the third point $(a,b_3)$ connected by involution to the two mentioned above
  is located in the inverse divisor $D_g^\ast = Z - D_{g+1} \geqslant 0$,
where $Z$ denotes the divisor of zeros of $\mathcal{G}$. 
\end{lem}
A proof is similar to the proof of Lemma~\ref{L:IZfactorize}.

\begin{rem}
The cases described in Lemmas~\ref{L:IZfactorize} and \ref{L:GZfactorize}
show that \eqref{Zf2g1} is not equivalent to \eqref{IfEqs1} when points of involution occur
in the divisor of zeros of an entire rational function $\mathcal{B}$.
In these cases $y$-coordinates of the points in involution are obtained from the curve equation.
\end{rem}

\subsection{Non-special divisors in terms of entire functions}
\begin{theo}\label{T:DgIG}
Suppose, the entire rational functions $\mathcal{I}$ and $\mathcal{G}$ of orders $2g$ and $2g+1$, respectively,  
vanish on the same degree $g$ strictly non-special  divisor $D_g$ on a trigonal curve \eqref{C3def}
of genus $g$, that is $D_g$ serves as the only solution of the system 
\begin{gather}\label{IGdefDg}
\mathcal{I}(x,y) = 0,\qquad \mathcal{G}(x,y) = 0.
\end{gather}
Then $D_g$ is equivalently defined by the system 
\begin{gather}\label{HIdefDg}
\mathcal{H}(x) = 0,\qquad \mathcal{I}(x,y) = 0,
\end{gather}
or
\begin{gather}\label{HGdefDg}
\mathcal{H}(x) = 0,\qquad \mathcal{G}(x,y) = 0,
\end{gather}
where $\mathcal{H}$ is a polynomial of degree $g$ such that
\begin{gather}\label{HDefC3}
\begin{aligned}
 \mathcal{H}(x) &=  \gamma_{y} (x) \mathcal{I}(x,y)  - \alpha_{y} (x) \mathcal{G}(x,y) \\
 &=   \gamma_{y} (x)  \alpha_{x} (x) -  \alpha_{y} (x)   \gamma_{x} (x).
 \end{aligned}
\end{gather}
$D_g$ is defined uniquely by any of the systems \eqref{IGdefDg}--\eqref{HGdefDg}.
\end{theo}
\begin{proof}
From  \eqref{IDef3} and  \eqref{GDef3} one easily finds that the degree of $\mathcal{H}$ equals $g$,
that is $3\mFr$ on a curve of the type $\mathfrak{C}^1$, and $3\mFr +1$ on a curve of the type $\mathfrak{C}^2$.
Thus, the equation $\mathcal{H}(x)=0$ gives $g$ values of the $x$-coordinate, which define $3g$ points on the curve,
namely: $\{(x_k,y_k^{(i)}), i=1,2,3\}_{k=1}^g$, where $y_k^{(1)}$, $y_k^{(2)}$, $y_k^{(3)}$ 
are solutions of $f(x_k,y)=0$. The functions  $\mathcal{I}$ and $\mathcal{G}$ are linear in $y$, and so they allow to
single out the unique point which belongs to $D_g$  from every group of three points $(x_k,y_k^{(i)})$, 
$i=1,2,3$, connected by involution. Evidently, $D_g$ is defined uniquely.
\end{proof}

\begin{theo}\label{T:g1FGCns}
A strictly non-special divisor $D_{g+1}$ on a trigonal curve  \eqref{C3def}
of genus $g$ is uniquely defined by the system 
\begin{gather}\label{g1FGns}
\mathcal{F}(x) = 0,\qquad \mathcal{G}(x,y)=0,
\end{gather} 
where $\mathcal{F}$ is a  polynomial  in $x$ of degree $g+1$,
and $\mathcal{G}$ is an entire rational function of weight $2g+1$, the both vanishing on $D_{g+1}$.
\end{theo}
\begin{proof}
is similar to the proof of Theorem~\ref{T:DgIG}.  The equation $\mathcal{F}(x)=0$ gives $g+1$ 
values of the $x$-coordinate, which define $g+1$ groups of points in involution on the curve.
Then $\mathcal{G}$ singles out the unique point which belongs to $D_{g+1}$  
from every group of three points in involution.
\end{proof}

\section{The process of reduction of a degree $g+1$ divisor}\label{s:Main}
An entire function $\mathcal{I}$ of weight $2g$ has one of the following forms:
 \begin{subequations}\label{ImonomDef}
\begin{align}\label{ImonomDef1}
&\mathfrak{C}^{1}:& &\mathcal{I}(x,y) = y \sum_{k=0}^{\mFr-1} \alpha_{3\mFr-1-3k} x^{k} 
 + \sum_{k=0}^{2\mFr} \alpha_{6\mFr-3k}  x^{k},&\\   \label{ImonomDef2}
&\mathfrak{C}^{2}:& &\mathcal{I}(x,y) =  y \sum_{k=0}^{\mFr} \alpha_{3\mFr-3k} x^{k} 
 + \sum_{k=0}^{2\mFr} \alpha_{6\mFr+2-3k}  x^{k};&
\end{align}
\end{subequations}
and an entire function $\mathcal{G}$ of weight $2g+1$, respectively:
\begin{subequations}\label{GmonomDef}
\begin{align}\label{GmonomDef1}
&\mathfrak{C}^{1}:& &\mathcal{G}(x,y) = y \sum_{k=0}^{\mFr} \gamma_{3\mFr-3k} x^{k} 
 +  \sum_{k=0}^{2\mFr} \gamma_{6\mFr+1-3k} x^{k},&\\   \label{GmonomDef2}
&\mathfrak{C}^{2}:& &\mathcal{G}(x,y) = y \sum_{k=0}^{\mFr} \gamma_{3\mFr+1-3k} x^{k} 
 +  \sum_{k=0}^{2\mFr+1} \gamma_{6\mFr+3-3k} x^{k}.&
\end{align}
\end{subequations}
In what follows, we always denote the coefficients of $\mathcal{I}$  by $\alpha_k$, and
the coefficients of $\mathcal{G}$  by~$\gamma_k$.

Below in this section we construct a solution of the reduction problem on a curve of the type  $\mathfrak{C}^1$. 
A similar solution for a curve of the type $\mathfrak{C}^2$ is presented in Appendix~\ref{A:Main35Family}.

\subsection{Minimal function defining a degree $g+1$ divisor}\label{ss:gg1Red}
Let a divisor $D_{g+1}$ be composed from a reduced strictly non-special divisor
$D_g$ and an additional point $(\bar{x},\bar{y})$
which does not coincide with any point from $D_g$.
The divisor $D_g$ is defined by a polynomial $\mathcal{H}$ in $x$ of degree $g$ 
and an entire function $\mathcal{I}$ of weight $2g$, the both vanishing on $D_g$. 
According to Theorem~\ref{T:g1FGCns},  $D_{g+1}$ is defined by the
polynomial $\mathcal{F}$ in $x$ of degree $g+1$ 
\begin{gather*}
\mathcal{F}(x) = \mathcal{H}(x) (x-\bar{x}),
\end{gather*} 
and an entire function $\mathcal{G}$ of weight $2g+1$. The latter 
is obtained from the given $\mathcal{H}$ and  $\mathcal{I}$ by the formula 
\begin{gather}\label{GfromI}
\mathcal{G}(x,y) =  \mathcal{I}(x,y)  \mathcal{M}(x) + \mathcal{H}(x) \mathcal{N}(x,y),
\end{gather}
where $\mathcal{I}$ has the form \eqref{ImonomDef1} and the following notation is used
\begin{align*}
& \mathcal{H}(x) = \sum_{k=0}^{3\mFr} h_{9\mFr-3k} x^k,\\
& \mathcal{M}(x) =  \sum_{i=0}^{3\mFr-1} a_{9\mFr-3-3i} x^i,\\
& \mathcal{N}(x,y) = y \sum_{i=0}^{\mFr-2} b_{3\mFr-4-3i} x^i + \sum_{i=0}^{2\mFr-1} b_{6\mFr-3-3i} x^i
\equiv y \nu_y(x) + \nu_x(x).
\end{align*}

The expression \eqref{GfromI} for $\mathcal{G}$ is an expansion in monomials \eqref{Monom1}.
Since $\mathcal{G}$ has weight of $2g+1=6\mFr+1$,
all coefficients of the monomials $y x^i$ with $4\mFr-2 \geqslant i \geqslant \mFr+1$
and $x^i$ with $5\mFr-1 \geqslant i \geqslant 2\mFr+1$  vanish, which leads to a system of $6\mFr-3$ linear equations 
\begin{gather}\label{BpEqs}
\mathfrak{B} \mathfrak{p}=0
\end{gather}
in $6\mFr-1$ unknown parameters 
$$\mathfrak{p} = (a_{0},\, a_3,\, \dots, a_{9\mFr-3},\, 
b_{2},\, \dots,\, b_{3\mFr-4},\, b_{0},\, \dots, \, b_{6\mFr-3} )^t.$$ 
The coefficient matrix has a block form, namely:
\begin{gather}\label{BlockM}
\mathfrak{B} =
 \begin{pmatrix} \mathfrak{G}_a & \mathfrak{H}_a \\ \mathfrak{G}_b & \mathfrak{H}_b \end{pmatrix}.
\end{gather}
The rows are divided into two parts: $(\mathfrak{G}_a$, $\mathfrak{H}_a)$ corresponding to 
the vanishing coefficients of the monomials
from $y x^{4\mFr-2}$ to $y x^{\mFr+1}$, and
$(\mathfrak{G}_b$, $\mathfrak{H}_b)$ corresponding to the vanishing coefficients 
of the monomials from $x^{5\mFr-1}$ to $x^{2\mFr+1}$.
The columns are divided into two parts corresponding to $\{a_{0}$, $a_3$, \ldots, $a_{9\mFr-3}\}$ and 
$\{ b_{2}$, \ldots, $b_{3\mFr-4}$, $b_{0}$, \ldots, $b_{6\mFr-3}\}$, respectively.

Block $\mathfrak{G}_a$ is of size $(3\mFr-2) \times 3\mFr$ 
with two zero columns on the right, the remaining $(3\mFr-2) \times (3\mFr-2)$ part is a lower triangular 
 $\mFr$-diagonal matrix, namely:
\begin{gather*}
\mathfrak{G}_a = 
\begin{pmatrix} \alpha_2 & 0 & \cdots & 0 & 0 & \cdots & 0 & 0 & 0\\
\alpha_5 &  \alpha_2 & \ddots & \vdots & \vdots &  & \vdots & \vdots & \vdots \\
\vdots & \ddots & \ddots & 0 & 0 & \cdots & 0 & 0 & 0 \\
\alpha_{3\mFr-1} & \ddots & \alpha_5 &  \alpha_2 & 0  & \cdots & 0 & 0 & 0 \\
0 & \alpha_{3\mFr-1}  & \ddots & \alpha_5 & \ddots &  \ddots & \vdots & \vdots & \vdots\\
\vdots & & \ddots & \ddots &\ddots &\alpha_2 & 0 & 0 & 0 \\
0  & \dots & 0 & \alpha_{3\mFr-1}  & \dots  & \alpha_5  & \alpha_2 & 0 & 0
 \end{pmatrix}.
\end{gather*}
Block $\mathfrak{G}_b$ of size $(3\mFr-1) \times 3 \mFr$ has one zero column on the right,
and the remaining $(3\mFr-1) \times (3\mFr-1)$ part  is lower triangular $(2\mFr+1)$-diagonal:
\begin{gather*}
\mathfrak{G}_b = 
\begin{pmatrix} \alpha_0 & 0 & \cdots & 0 & 0 & \cdots & 0 & 0\\
\alpha_3 &  \alpha_0 & \ddots & \vdots & \vdots &  & \vdots & \vdots \\
\vdots & \ddots & \ddots & 0 & 0 & \cdots &0 & 0 \\
\alpha_{6\mFr} & \ddots & \alpha_3 &  \alpha_0 & 0  & \cdots & 0 & 0 \\
0 &  \alpha_{6\mFr} & \ddots & \alpha_3 & \ddots &  \ddots & \vdots & \vdots\\
\vdots & & \ddots & \ddots &\ddots &\alpha_0 & 0 & 0 \\
0  & \dots & 0 & \alpha_{6\mFr}  & \dots  & \alpha_3  & \alpha_0 & 0
 \end{pmatrix}.
\end{gather*}
Block $\mathfrak{H}_a$ of size $(3\mFr-2)\times (3\mFr-1)$ has $2\mFr$  zero columns on the right, 
and is of the form:
\begin{gather*}
\mathfrak{H}_a =  \begin{pmatrix} 
h_0 & 0 & \cdots & 0 & 0 & \cdots & 0 \\
h_3 & h_0 & \ddots & \vdots & 0 & \cdots & 0 \\
\vdots & \ddots & \ddots & 0 & \vdots & & \vdots \\
h_{3\mFr-6} & h_{3\mFr-9} & \ddots & h_0 & 0 & \cdots & 0 \\
h_{3\mFr-3} & h_{3\mFr-6}  & \cdots & h_3 & 0 &  \cdots & 0 \\
\vdots & \vdots & & \vdots & \vdots & & \vdots \\
h_{9\mFr-9} & h_{9\mFr-12} & \cdots & h_{6\mFr-3} & 0 & \cdots &0
\end{pmatrix}.
\end{gather*}
Block $\mathfrak{H}_b$ of size $(3\mFr-1) \times (3\mFr-1)$ has $\mFr-1$ zero columns on the left,
and is of the form:
\begin{gather*}
\mathfrak{H}_b =  \begin{pmatrix} 
0 & \cdots & 0 & h_0 & 0 & \cdots & 0 \\
0 & \cdots & 0 & h_3 & h_0 & \ddots & \vdots \\
\vdots & &  \vdots & \vdots & \vdots & \ddots & 0 \\
0 & \cdots & 0 & h_{6\mFr-3} & h_{6\mFr-6} & \cdots & h_0 \\
\vdots & &  \vdots & \vdots & \vdots &  & \vdots \\
0 & \cdots & 0 & h_{9\mFr-6} & h_{9\mFr-9} & \cdots & h_{3\mFr-3}
\end{pmatrix},
\end{gather*}
Evidently, the kernel of $\mathfrak{B}$  is two-dimensional, and $a_{9\mFr-3}$ 
is a free parameter since the $3\mFr$-th column of
blocks $\mathfrak{G}_a$ and $\mathfrak{G}_b$ is zero.
In addition to an arbitrary $a_{9\mFr-3}$, the parameters $\mathfrak{p}$
are defined up to a common constant multiple, which is~$a_0$.

We introduce the function $\widetilde{\mathcal{G}}(x,y)$ defined by \eqref{GfromI}
with $a_{9\mFr-3}=0$, and assign
\begin{gather*}
a_{9\mFr-3} = -\frac{\widetilde{\mathcal{G}}(\bar{x},\bar{y})}{\mathcal{I}(\bar{x},\bar{y})}.
\end{gather*}
Recall that $(\bar{x},\bar{y})$ does not coincide with any point from $D_g$. Thus, 
$\mathcal{I}(\bar{x},\bar{y}) \neq 0$, and so $a_{9\mFr-3}$ is well-defined.
Note, that $a_0$ arises as an arbitrary factor of $\widetilde{\mathcal{G}}(x,y)$.
It is convenient to use it to rid off the denominators in all coefficients.
Then  the desired function $\mathcal{G}$
vanishing on $D_{g+1} = D_g + (\bar{x},\bar{y})$ is given by the formula
\begin{gather}\label{Gexpr}
\mathcal{G}(x,y) = \widetilde{\mathcal{G}}(x,y) 
- \frac{\widetilde{\mathcal{G}}(\bar{x},\bar{y})}{\mathcal{I}(\bar{x},\bar{y})} \mathcal{I}(x,y).
\end{gather}

\begin{rem}\label{R:SamePs}
If $(\bar{x},\bar{y})$ coincides with one of the points from $D_g$, then  $\mathcal{I}(\bar{x},\bar{y}) = 0$.
It means $\widetilde{\mathcal{G}}(\bar{x},\bar{y}) / \mathcal{I}(\bar{x},\bar{y})$ should be replaced with
$\lim_{(x,y) \to (\bar{x},\bar{y})} \big( \widetilde{\mathcal{G}}(x,y) / \mathcal{I}(x,y) \big)$.
The latter exists due to $(\bar{x},\bar{y})$ is a simple zero of $\mathcal{I}$ and $\widetilde{\mathcal{G}}$, though 
this situation requires a special attention in the process of implementation.
Considering an example in section~\ref{s:C34impl}, we show how to avoid such a difficulty.
\end{rem}

\begin{exam}
In the case of a $C_{3,4}$, we have
\begin{gather*}
\mathfrak{B} = \begin{pmatrix} 
 \alpha_2 & 0 & 0 & 0 & 0 \\
 \alpha_0 & 0 & 0 & h_0 & 0 \\
 \alpha_3 & \alpha_0 & 0 & h_3 & h_0 
 \end{pmatrix},\qquad
\mathfrak{p} = 
\big(a_0,\, a_3,\,  a_6,\, b_0,\, b_3\big)^t,
\end{gather*}
and so
\begin{multline}\label{GfromIHC34}
 \widetilde{\mathcal{G}}(x,y) = \alpha_2^{-1} \big(h_0 x \mathcal{I}(x,y) - \alpha_0 \mathcal{H}(x)\big)\\
= \alpha_2^{-1} \big(h_0 \alpha_2  y x + (h_0 \alpha_3 - h_3 \alpha_3) x^2 +
(h_0 \alpha_6 - h_6 \alpha_0) x - h_9 \alpha_0\big).
\end{multline}
\end{exam}

\begin{rem}
Composing the expression \eqref{GfromI} for $\mathcal{G}$, 
we multiply $\mathcal{I}$ by a 
polynomial $\mathcal{M}$ of degree equal to $\deg (\mathcal{H}) -1=g-1$,
and multiply  $\mathcal{H}$ by an entire function $\mathcal{N}$ of weight equal to 
$\wgt (\mathcal{I}) - \wgt (x)$, that is 
$\deg (\nu_y) = \deg (\alpha_y) -1$ and $\deg (\nu_x) = \deg (\alpha_x) -1$.
\end{rem}

\subsection{Inverse to a degree $g+1$ divisor}\label{ss:g1Div}
Let $D_{g+1}$ be a degree $g+1$ strictly non-special divisor
which is defined by a polynomial $\mathcal{F}$ of degree $g+1$ and
an entire function $\mathcal{G}$ of weight $2g+1$ of the form \eqref{GmonomDef1}, the both
vanishing on $D_{g+1}$.
In this subsection we find the divisor $D^\ast_g$ inverse to $D_{g+1}$.

It is sufficient to define  $D^\ast_g$ by the given function $\mathcal{G}$ and
a polynomial $\mathcal{H}^\ast$ in $x$ of degree~$g$, as stated in Theorem~\ref{T:DgIG}.
The polynomial $\mathcal{H}^\ast$ is obtained by the formula
\begin{gather}\label{HDef}
\mathcal{H}^\ast(x) = \frac{\mathcal{Z}_{2g+1}(x)}{\mathcal{F}(x)} 
= \frac{\gamma_y(x)^3 f\big(x,-\gamma_x(x)/\gamma_y(x)\big)}{\mathcal{F}(x)},
\end{gather}
where $\mathcal{Z}_{2g+1}$, defined by \eqref{Zgp}, is divisible by $\mathcal{F}$
due to  $D_{g+1}$ is a part of the divisor of zeros of  $\mathcal{G}$.

Therefore, $D^\ast_g$ is defined by
\begin{gather*}
\mathcal{H}^\ast(x) = 0,\qquad \mathcal{G}(x,y)=0.
\end{gather*}

\subsection{Minimal function defining a degree $g$ divisor}\label{ss:g1gRed}
Let $D^\ast_g $ be a degree $g$ strictly non-special divisor 
defined by a polynomial $\mathcal{H}^\ast$ of degree $g$ and an entire
function $\mathcal{G}$ of weight $2g+1$,  the both vanishing on $D^\ast_g$. Now we
find the minimal function $\mathcal{I}$ defining $D^\ast_g$.
In other words, we reduce the given function $\mathcal{G}$ of weight $2g+1$ 
to a function $\mathcal{I}$ of weight~$2g$, taking into account that $D^\ast_g$
is a part of the divisor of zeros of~$\mathcal{H}^\ast$.

We construct the required function $\mathcal{I}$ by the formula
\begin{gather}\label{IfromG}
\mathcal{I}(x,y) =  \mathcal{G}(x,y)  \mathcal{M}(x) + \mathcal{H}^\ast(x) \mathcal{N}(x,y),
\end{gather}
where $\mathcal{G}$ has the form \eqref{GmonomDef1}, and
\begin{align*}
&\mathcal{H}^\ast (x) = \sum_{k=0}^{3\mFr} \bar{h}_{9\mFr-3k} x^k, \\
& \mathcal{M}(x) =  \sum_{i=0}^{3\mFr-1} a_{9\mFr-3-3i} x^i,\\
& \mathcal{N}(x,y) = y \sum_{i=0}^{\mFr-1} b_{3\mFr-3-3i} x^i + \sum_{i=0}^{2\mFr-1} b_{6\mFr-2-3i} x^i
\equiv y \nu_y(x) + \nu_x(x).
\end{align*}
The expression for $\mathcal{I}$ contains monomials \eqref{Monom1} up to the weight $2g=6\mFr$. Thus,
all coefficients of the monomials $y x^i$ with $4\mFr-1 \geqslant i \geqslant \mFr$
and $x^i$ with $5\mFr-1 \geqslant i \geqslant 2 \mFr+1$  vanish, which leads to a system of $6\mFr-1$ 
linear equations of the form \eqref{BpEqs}
in $6\mFr$ unknown parameters 
$$\mathfrak{p} = (a_{0},\, a_3,\, \dots, a_{9\mFr	-3},\, 
b_{0},\, b_3,\, \dots,\, b_{3\mFr-3},\, b_{1},\, b_4,\, \dots, \, b_{6\mFr-2} )^t.$$ 
The coefficient matrix has the block form \eqref{BlockM}.
Row $(\mathfrak{G}_a$, $\mathfrak{H}_a)$ corresponds to 
the vanishing coefficients of the monomials
from $y x^{4\mFr-1}$ to $y x^{\mFr}$; and row
$(\mathfrak{G}_b$, $\mathfrak{H}_b)$ corresponds to 
the vanishing coefficients of the monomials from $x^{5\mFr-1}$ to $x^{2\mFr+1}$.
Columns are divided into two parts corresponding to $\{a_{0}$, $a_3$, \ldots, $a_{9\mFr-3}\}$ and 
$\{b_{0}$,  \ldots, $b_{3\mFr-3}$, $b_{1}$, \ldots, $b_{6\mFr-2}\}$.

Block $\mathfrak{G}_a$ of size $3\mFr \times 3\mFr$  is a lower triangular $(\mFr+1)$-diagonal matrix:
\begin{gather*}
\mathfrak{G}_a = 
 \begin{pmatrix} \gamma_0 & 0 & \cdots & 0 & 0 & \cdots & 0\\
\gamma_3 &  \gamma_0 & \ddots & \vdots & \vdots &  & \vdots \\
\vdots & \ddots & \ddots & 0  & 0 & \cdots & 0 \\
\gamma_{3\mFr}  &  \ddots & \gamma_3 &  \gamma_0 & 0 &  \cdots & 0\\
0 & \gamma_{3\mFr} & \ddots &  \gamma_3 & \ddots & \ddots & \vdots\\
\vdots & & \ddots & \ddots & \ddots &  \gamma_0 &  0 \\
0 & \cdots & 0 & \gamma_{3\mFr} &  \cdots & \gamma_3 &  \gamma_0
 \end{pmatrix}.
\end{gather*}
Block $\mathfrak{G}_b$ is of size $(3\mFr-1) \times 3\mFr$,
and has one zero column on the right; the remaining $(3\mFr-1) \times (3\mFr-1)$
matrix is lower triangular $(2\mFr+1)$-diagonal: 
\begin{gather*}
\mathfrak{G}_b = 
  \begin{pmatrix} \gamma_1 & 0 & \cdots & 0 & 0 & \cdots & 0 & 0 \\
\gamma_4 &  \gamma_1 & \ddots & \vdots & \vdots &  & \vdots & \vdots \\
\vdots & \ddots & \ddots & 0 & 0 & \cdots & 0 & 0 \\
\gamma_{6\mFr+1}  &  \ddots & \gamma_4 &  \gamma_1 & 0 &  \cdots & 0 & 0\\
0 & \gamma_{6\mFr+1} & \ddots & \gamma_4 & \ddots & \ddots & \vdots & \vdots\\
\vdots & & \ddots & \ddots & \ddots &   \gamma_1 & 0 & 0  \\
0 & \cdots & 0 & \gamma_{6\mFr+1} &  \cdots & \gamma_4 &  \gamma_1 & 0
 \end{pmatrix}.
\end{gather*}
Block $\mathfrak{H}_a$ of size $3\mFr\times 3\mFr$ has $2\mFr$ zero columns on the right,
and block $\mathfrak{H}_b$ of size $(3\mFr-1)\times 3\mFr$ has $\mFr$ zero columns on the left, namely:
\begin{gather*}
\mathfrak{H}_a =  \begin{pmatrix} 
\bar{h}_0 & 0 & \cdots & 0 & 0 & \cdots & 0 \\
\bar{h}_3 & \bar{h}_0 & \ddots & \vdots & 0 & \cdots & 0 \\
\vdots & \ddots & \ddots & 0 & \vdots & & \vdots \\
\bar{h}_{3\mFr-3} & \bar{h}_{3\mFr-6} & \ddots & \bar{h}_0 & 0 & \cdots & 0 \\
\bar{h}_{3\mFr} & \bar{h}_{3\mFr-3}  & \cdots & \bar{h}_3 & 0 &  \cdots & 0 \\
\vdots & \vdots & & \vdots & \vdots & & \vdots \\
\bar{h}_{9\mFr-3} & \bar{h}_{9\mFr-6} & \cdots & \bar{h}_{6\mFr} & 0 & \cdots &0
\end{pmatrix},\\
\mathfrak{H}_b =  \begin{pmatrix} 
0 & \cdots & 0 & \bar{h}_0 & 0 & \cdots & 0 \\
0 & \cdots & 0 & \bar{h}_3 & \bar{h}_0 & \ddots & \vdots\\
\vdots & &  \vdots & \vdots & \vdots & \ddots & 0 \\
0 & \cdots & 0 & \bar{h}_{6\mFr-3} & \bar{h}_{6\mFr-6} & \cdots & \bar{h}_0\\
\vdots & &  \vdots & \vdots & \vdots &  & \vdots \\
0 & \cdots & 0 & \bar{h}_{9\mFr-6} & \bar{h}_{9\mFr-9} & \cdots & \bar{h}_{3\mFr-3}
\end{pmatrix}.
\end{gather*}
The kernel of $\mathfrak{B}$  is one-dimensional, so the 
parameters $\mathfrak{p}$
are defined up to a constant multiple, which is again $a_0$. 
We suggest to rid off denominators of all coefficients of $\mathcal{I}$
by means of $a_0$.

\begin{exam}
In the case of $C_{3,4}$, we have
\begin{gather*}
\mathfrak{B} = \begin{pmatrix} 
 \gamma_0 & 0 & 0 & \bar{h}_0 & 0 & 0\\
 \gamma_3 & \gamma_0 & 0 & \bar{h}_3 & 0 & 0 \\
 0 & \gamma_3 & \gamma_0 & \bar{h}_6 & 0 & 0 \\
 \gamma_1 & 0 & 0 & 0 & \bar{h}_0 & 0 \\
 \gamma_4 & \gamma_1 & 0  & 0 & \bar{h}_3 & 0
\end{pmatrix},\qquad
\mathfrak{p} = 
\big(a_0,\, a_3,\,  a_6,\, b_0,\, b_1,\, b_4\big)^t,
\end{gather*}
and \eqref{IfromG}  gets the form
\begin{multline}\label{IfromGHC34r}
\sfit{C}_7 \mathcal{I}(x,y) =  \Big(\bar{h}_0 \gamma_0^2 x^2 + (\bar{h}_3\gamma_0^2 - \bar{h}_0 \gamma_0 \gamma_3) x 
+ \bar{h}_6  \gamma_0^2 - \bar{h}_3 \gamma_0 \gamma_3 + \bar{h}_0 \gamma_3^2  \Big) \mathcal{G}(x,y)  \\
- \big(\gamma_0^3 y + \gamma_0^2 \gamma_1 x + \gamma_0^2 \gamma_4 - \gamma_0 \gamma_3 \gamma_1  \big) 
\mathcal{H}(x),
\end{multline}
where $\sfit{C}_7$ is a normalisation constant, which has S\={a}to weight equal to $7$,
as seen from the equality. Substituting $\mathcal{G}$ and $\mathcal{H}$, we find
\begin{subequations}\label{IfromGHC34s}
\begin{multline}
\sfit{C}_7 \mathcal{I}(x,y) =  \bar{h}_0 (\gamma_1 \gamma_3^2 
- \gamma_0 \gamma_3 \gamma_4 + \gamma_0^2 \gamma_7) x^2
+ y \big(\bar{h}_0 \gamma_3^3 - \bar{h}_3 \gamma_0 \gamma_3^2 
+ \bar{h}_6 \gamma_0^2 \gamma_3 \\ - \bar{h}_9 \gamma_0^3 \big)  
+ x \big(\bar{h}_0 (\gamma_3^2 \gamma_4 - \gamma_0 \gamma_3 \gamma_7)
- \bar{h}_3 (\gamma_0 \gamma_3 \gamma_4 - \gamma_0^2 \gamma_7)
+ \bar{h}_6 \gamma_0 \gamma_1 \gamma_3 - \bar{h}_9 \gamma_0^2 \gamma_1\big) \\
+ \bar{h}_0 \gamma_3^2 \gamma_7 - \bar{h}_3 \gamma_0\gamma_3 \gamma_7
+ \bar{h}_6 \gamma_0^2 \gamma_7 + \bar{h}_9 (\gamma_0\gamma_1\gamma_3 - \gamma_0^2 \gamma_4),
\end{multline}
where
\begin{gather}
\sfit{C}_7 = \gamma_1 \gamma_3^2 - \gamma_0 \gamma_3 \gamma_4 + \gamma_0^2 \gamma_7.
\end{gather}
\end{subequations}
\end{exam}

\begin{rem}
Composing the expression \eqref{IfromG} for $\mathcal{I}$, we multiply $\mathcal{G}$ by 
a polynomial $\mathcal{M}$ of degree $\deg (\mathcal{H}) -1=g-1$, and multiply
$\mathcal{H}$ by an entire function $\mathcal{N}$ of weight equal to 
$\wgt \mathcal{G} - \wgt x$, that is 
$\deg \nu_y = \deg \gamma_y -1$ and $\deg \nu_x = \deg \gamma_x -1$.
\end{rem}

\subsection{Inverse to a degree $g$ divisor}\label{ss:gInv}
Let $D^\ast_g$ be a degree $g$ strictly non-special divisor defined by
a polynomial $\mathcal{H}^\ast$ in $x$ of degree $g$
and an entire function $\mathcal{I}$ of weight $2g$, the both vanishing on $D^\ast_g$.
We find the divisor $\widetilde{D}_g$ inverse to $D^\ast_g$.

According to Theorem~\ref{T:DgIG}, $\widetilde{D}_g$ as a degree $g$ divisor
is defined by the given function~$\mathcal{I}$,
and a polynomial $\widetilde{\mathcal{H}}$ in $x$ of degree~$g$. 
We obtain $\widetilde{\mathcal{H}}$ by the formula:
\begin{gather}\label{HinvDef}
\widetilde{\mathcal{H}}(x) = \frac{\mathcal{Z}_{2g}(x)}{\mathcal{H}^\ast(x)} 
= \frac{\alpha_y(x)^3 f\big(x,-\alpha_x(x)/\alpha_y(x)\big)}{\mathcal{H}^\ast(x)},
\end{gather}
where $\mathcal{Z}_{2g}$, defined by \eqref{Zgp}, is divisible by $\mathcal{H}^\ast$
due to  $D_{g}^\ast$ is a part of the divisor of zeros of  $\mathcal{I}$.

Therefore, $\widetilde{D}_g$ is defined by
\begin{gather*}
\widetilde{\mathcal{H}}(x) = 0,\qquad \mathcal{I}(x,y)=0.
\end{gather*}

\section{Reduction algorithm}\label{s:RedAlg}
In this section we present an iterative algorithm which solves the reduction problem for a
degree $g+p$ divisor $D_{g+p}$ with $p>0$. We suppose that $D_{g+p}$ is non-special and contains
no groups of three points in involution.

We start with an arbitrary degree $g$ divisor $D_{g}\subset D_{g+p}$,
that is $D_{g} = \sum_{k=1}^{g} (x_k,\, y_k)$.
We define $D_{g}$ by the system
\begin{gather*}
\mathcal{H}(x) = 0,\qquad  \mathcal{I}(x,y) = 0.
\end{gather*} 
The polynomial $\mathcal{H}$ in $x$ of degree $g$ has zeros $\{x_k\}_{k=1}^{g}$, and
the leading coefficient equals $1$.
The entire function $\mathcal{I}$ of weight $2g$ is defined by \eqref{R2gMatr}.

We denote by $D_{\text{res}}$ the remaining part of $D_{g+p}$,
that is $D_{\text{res}} = D_{g+p} -  D_{g} \geqslant 0$.
\begin{enumerate}
\renewcommand{\labelenumi}{\textbf{Step \arabic{enumi}}}
\item Let  $(\bar{x},\bar{y})$ be an arbitrary point  from $D_{\text{res}}$, not coinciding with
any point from $D_g$.
Construct an entire function $\mathcal{G}$ of weight $2g+1$ with the help of \eqref{GfromI}, 
see subsection~\ref{ss:gg1Red} for more detail. Assign $\mathcal{F}(x) = \mathcal{H}(x) (x-\bar{x})$.
Then a new divisor $D_{g+1}$ is defined by the system
\begin{gather*}
\mathcal{F}(x) = 0,\qquad  \mathcal{G}(x,y)=0.
\end{gather*}

\item 
Find a degree $g$ polynomial $\mathcal{H}^\ast$ by \eqref{HDef},
 see subsection~\ref{ss:g1Div}. Then the system
\begin{gather*}
\mathcal{H}^\ast (x) = 0,\qquad \mathcal{G}(x,y)=0.
\end{gather*}
defines the divisor $D^\ast_g$  complement to $D_{g+1}$
in the divisor of zeros of $\mathcal{G}$. 

\item  By \eqref{IfromG}  construct an entire rational function $\mathcal{I}^\ast$ of weight $2g$, 
which is the minimal function 
defining  $D^\ast_g$ from the previous step. A new definition of $D^\ast_g$ is
\begin{gather*}
\mathcal{H}^\ast(x) = 0,\qquad  \mathcal{I}^\ast(x,y) = 0.
\end{gather*}
The procedure of reducing $\mathcal{G}$ to $\mathcal{I}^\ast$ 
is explained in subsection~\ref{ss:g1gRed}.

\item 
Find a degree $g$ polynomial $\widetilde{\mathcal{H}}$ by \eqref{HinvDef},
 see subsection~\ref{ss:gInv}. 
The system 
\begin{gather*}
\widetilde{\mathcal{H}}(x) = 0,\qquad \mathcal{I}^\ast(x,y)=0
\end{gather*}
defines the reduced divisor $\widetilde{D}_g$ equivalent
to $D_{g+1}$ from Step~1.
\end{enumerate}
These new functions $\widetilde{\mathcal{H}}$ and $\mathcal{I}^\ast$ 
are initial for the next iteration.

\medskip
\textbf{Summary.} In the first iteration of the algorithm we construct a degree $g+1$ divisor $D_{g+1}$ by adding 
one point from  $D_{\text{res}}$ to $D_g$,
then reduce $D_{g+1}$  to the equivalent degree $g$ divisor $\widetilde{D}_g$ (a reduced divisor).
In each further iteration we construct  a new degree $g+1$ divisor $D_{g+1}$ by 
adding one of the remaining points from  $D_{\text{res}}$ to the reduced divisor
$\widetilde{D}_g$ obtained in the previous iteration, then find a new reduced divisor
equivalent to the new degree $g+1$ divisor $D_{g+1}$.

\medskip
Addition of two degree $g$ non-special  divisors is realized as follows. We start with one of the two divisors,
and add another divisor point by point.
In the case of  the standard cryptography problem:  a degree $g$ divisor $D_g$ is added to itself.
This causes some computational issue in Step 1,
since $\mathcal{I}(\bar{x},\bar{y}) = 0$, see Remark~\ref{R:SamePs}. 
We avoid such an issue by starting the first iteration from Step 2.
The polynomial $\mathcal{F}$ of degree $g+1$ and the entire function $\mathcal{G}$ of weight $2g+1$
are constructed directly from a degree $g+1$ divisor $D_{g+1}$, which contains at least two equal points.
Remark~\ref{R:repPs} explains how to construct $\mathcal{G}$ in the case of repeated points.
Such a situation occurring in the next iterations is practically impossible.

\section{Implementation of reduction algorithm on $C_{3,4}$}\label{s:C34impl}
In this section the reduction algorithm on $C_{3,4}$ is presented explicitly, 
that is in the form ready for implementation. 
Polynomial division in Steps 2 and 4 is performed in terms of coefficients, 
as explained in Appendix~\ref{A:PolyDiv}.
Note, that the indices of coefficients display their S\={a}to wights. All expressions below
respect the S\={a}to wights, that is all summands in any expression have equal weights,
and the weights on the both sides of relations are equal. This helps to verify the correctness of relations.

We start each iteration, except possibly the first, with 
\begin{align*}
&\mathcal{H}(x) = h_0 x^3 + h_{3} x^2 + h_{6} x + h_{9},\\
&\mathcal{I}(x,y) = \alpha_0 x^2 + \alpha_2 y + \alpha_3 x + \alpha_{6},
\end{align*}
which define the reduced divisor $D_3$  from the previous iteration.
We denote the point which is added by $(\bar{x},\,\bar{y})$.
\begin{enumerate}
\item[\textbf{Step 1.}]
Find
\begin{align*}
&\mathcal{F}(x) = f_0 x^4 + f_{3} x^3 + f_{6} x^2 + f_{9} x + f_{12},\\
&\widetilde{\mathcal{G}}(x,y) = \bar{\gamma}_0 x y  + \bar{\gamma}_1 x^2 
+ \bar{\gamma}_3 y + \bar{\gamma}_4 x +\bar{\gamma}_7,
\end{align*}
where $\mathcal{F}(x) = (x-\bar{x})\mathcal{H}(x)$, which implies
\begin{equation*}
f_{0} = h_0,\qquad f_{3i+3} = h_{3i+3} - h_{3i} \bar{x},\ i=0,\,\dots,\,2, \qquad
f_{12} = - h_{9} \bar{x},
\end{equation*}
and
\begin{align*}
&\bar{\gamma}_0 = h_0,&
&\bar{\gamma}_1 = \alpha_2^{-1}  (h_0 \alpha_3 - h_3 \alpha_0)\\
&\bar{\gamma}_3 = 0 &
&\bar{\gamma}_4 =  \alpha_2^{-1} (h_0 \alpha_6 - h_6 \alpha_0)\\
&& &\bar{\gamma}_7 =-  \alpha_2^{-1} h_9 \alpha_0.
\end{align*}
Let $c = \widetilde{\mathcal{G}}(\bar{x},\bar{y})/ \mathcal{I}(\bar{x},\bar{y})$. Then
\begin{align*}
&\mathcal{G}(x,y) = \gamma_0 x y  + \gamma_1 x^2 
+ \gamma_3 y + \gamma_4 x + \gamma_7,
\end{align*}
where
\begin{align*}
&\gamma_0 = \bar{\gamma}_0, \qquad
\gamma_{n} = \bar{\gamma}_{n} - c \alpha_{n-1},\quad n = 1,\, 3,\, 4,\, 7.
\end{align*}

\item[\textbf{Step 2.}] 
Find
\begin{equation*}
\mathcal{H}^\ast (x) = \bar{h}_0 x^3 + \bar{h}_{3} x^2 + \bar{h}_{6} x + \bar{h}_{9},
\end{equation*}
with
\begin{gather*}
\bar{h}_{9-3 k} = \sum_{i=0}^{k} (-1)^{i} d_{-12-3i} c_{21-3(k-i)},\qquad k=0,\dots,3,
\end{gather*}
where $d_{-12-3i}$ produced by \eqref{dExpr} have the form
\begin{align*}
&d_{-12} =  f_{12}^{-1},\\
&d_{-15} = f_{12}^{-2} f_9,\\
&d_{-18} = f_{12}^{-3} f_9^2 - f_{12}^{-2}  f_6,\\
&d_{-21} = f_{12}^{-4} f_9^3 - 2 f_{12}^{-3} f_6 f_9 + f_{12}^{-2} f_3,
\end{align*}
and $c_{12}$, $c_{15}$, $c_{18}$, $c_{21}$ are coefficients of $\mathcal{Z}_{7}$ defined by \eqref{Zgp},
namely
\begin{align*}
&c_{21} = - \gamma_7^3 + \gamma_3 \gamma_7^2 \lambda_4  - \gamma_3^2 \gamma_7  \lambda_8 
+ \gamma_3^3 \lambda_{12},\\
&c_{18} =  - 3 \gamma_4 \gamma_7^2 + \gamma_3 \gamma_7^2 \lambda_1   
+ (2\gamma_3 \gamma_4 \gamma_7 + \gamma_0 \gamma_7^2) \lambda_4 
-  \gamma_3^2 \gamma_7  \lambda_5 \\
&\quad \quad
- (\gamma_3^2 \gamma_4 + 2 \gamma_0 \gamma_3 \gamma_7) \lambda_8
+ \gamma_3^3 \lambda_9 + 3 \gamma_0 \gamma_3^2 \lambda_{12},\\
&c_{15} = -3 \gamma_4^2 \gamma_7 - 3 \gamma_1 \gamma_7^2
+ (2\gamma_3 \gamma_4 \gamma_7 + \gamma_0 \gamma_7^2) \lambda_1 
- \gamma_3^2 \gamma_7 \lambda_2  \\
&\quad \quad
+ (\gamma_3 \gamma_4^2 + 2 \gamma_1 \gamma_3 \gamma_7 
+ 2 \gamma_0 \gamma_4 \gamma_7) \lambda_4
- (\gamma_3^2 \gamma_4 + 2 \gamma_0 \gamma_3 \gamma_7) \lambda_5   + \gamma_3^3 \lambda_6\\
&\quad \quad
- (\gamma_1 \gamma_3^2 + 2 \gamma_0 \gamma_3 \gamma_4 
+ \gamma_0^2\gamma_7) \lambda_8 + 3 \gamma_0 \gamma_3^2  \lambda_9 
+ 3 \gamma_0 \gamma_3 \lambda_{12},\\
&c_{12} =  - \gamma_4^3  - 6 \gamma_1 \gamma_4 \gamma_7
+ (\gamma_3 \gamma_4^2 + 2 \gamma_1 \gamma_3 \gamma_7 
+ 2 \gamma_0 \gamma_4 \gamma_7) \lambda_1  \\
&\quad \quad - (\gamma_3^2 \gamma_4 + 2 \gamma_0 \gamma_3 \gamma_7) \lambda_2 + \gamma_3^3 \lambda_3 
+ (2\gamma_1 \gamma_3 \gamma_4 + \gamma_0 \gamma_4^2 + 2 \gamma_0 \gamma_1 \gamma_7) \lambda_4\\
&\quad \quad 
- (\gamma_1 \gamma_3^2 + 2 \gamma_0 \gamma_3 \gamma_4 
+ \gamma_0^2\gamma_7) \lambda_5 + 3 \gamma_0 \gamma_3^2 \lambda_6\\
&\quad\quad  - (2\gamma_0\gamma_1 \gamma_3 + \gamma_0^2 \gamma_4) \lambda_8 
+ 3 \gamma_0 \gamma_3 \lambda_9 + \gamma_0^3\lambda_{12}.
\end{align*}

\item[\textbf{Step 3.}] 
Find the new
\begin{equation*}
\mathcal{I}^\ast(x,y) = \bar{\alpha}_0 x^2 + \bar{\alpha}_2 y + \bar{\alpha}_3 x + \bar{\alpha}_{6},
\end{equation*}
where
\begin{align*}
&\bar{\alpha}_0 = \bar{h}_0,\\
&\bar{\alpha}_2 = - \sfit{C}_7^{-1}\big(\gamma_0 ( \gamma_0^2 \bar{h}_9 
- \gamma_0 \gamma_3 \bar{h}_6 + \gamma_3^2 \bar{h}_3) 
- \gamma_3^2 \gamma_3 \bar{h}_0 \big), \\
&\bar{\alpha}_3 = - \sfit{C}_7^{-1}\big(\gamma_1( \gamma_0^2 \bar{h}_9 - \gamma_0 \gamma_3 \bar{h}_6) 
-  ( \gamma_0^2 \gamma_7 - \gamma_0 \gamma_3 \gamma_4) \bar{h}_3 
 -  (- \gamma_0 \gamma_3\gamma_7 + \gamma_3^2 \gamma_4) \bar{h}_0 \big), \\
&\bar{\alpha}_6 = - \sfit{C}_7^{-1}\big( ( \gamma_0^2 \gamma_4 - \gamma_0 \gamma_3 \gamma_1) \bar{h}_9
- \gamma_7 ( \gamma_0^2 \bar{h}_6 - \gamma_0 \gamma_3 \bar{h}_3 + \gamma_3^2 \bar{h}_0)\big),\\
&\sfit{C}_7 = \gamma_1\gamma_3^2 + \gamma_0 \gamma_3 \gamma_4 + \gamma_0^2 \gamma_7.
\end{align*}

\item[\textbf{Step 4.}] 
Find the new
\begin{equation*}
\widetilde{\mathcal{H}}(x) = \tilde{h}_0 x^3 + \tilde{h}_{3} x^2 + \tilde{h}_{6} x + \tilde{h}_{9},
\end{equation*}
with
\begin{gather*}
\tilde{h}_{9-3 k} = \sum_{i=0}^{k} (-1)^{i} d_{-9-3i} c_{18-3(k-i)},\qquad k=0,\dots,3, 
\intertext{where}
\begin{split} \notag
&d_{-9} = \bar{h}_{9}^{-1},\\
&d_{-9-3j} = \sum_{i=1}^{j} (-1)^{i-1} \bar{h}_{9}^{-1} \bar{h}_{9-3i} d_{-9-3(j-i)},\qquad j = 1,\dots, 3,
\end{split}
\end{gather*}
and 
$c_{9}$, $c_{12}$, $c_{15}$, $c_{18}$ are coefficients of $\mathcal{Z}_{6}$ defined by \eqref{Zgp},
namely
\begin{align*}
&c_{18} = - \bar{\alpha}_6^3 + \bar{\alpha}_2 \bar{\alpha}_6^2 \lambda_4 - \bar{\alpha}_2^2 \bar{\alpha}_6 \lambda_8 + \bar{\alpha}_2^3 \lambda_{12},\\
&c_{15} = - 3 \bar{\alpha}_3 \bar{\alpha}_6^2 + \bar{\alpha}_2 \bar{\alpha}_6^2 \lambda_1 
+ 2 \bar{\alpha}_2 \bar{\alpha}_3 \bar{\alpha}_6 \lambda_4 
- \bar{\alpha}_2^2 \bar{\alpha}_6 \lambda_5 - \bar{\alpha}_2^2 \bar{\alpha}_3 \lambda_8 + \bar{\alpha}_2^3 \lambda_9,\\
&c_{12} = - 3 \bar{\alpha}_3^2 \bar{\alpha}_6 - 3 \bar{\alpha}_0 \bar{\alpha}_6^2
+ 2 \bar{\alpha}_2 \bar{\alpha}_3 \bar{\alpha}_6 \lambda_1 - \bar{\alpha}_2^2 \bar{\alpha}_6 \lambda_2 \\
&\quad \quad + ( \bar{\alpha}_2 \bar{\alpha}_3^2 + 2 \bar{\alpha}_0 \bar{\alpha}_2 \bar{\alpha}_6) \lambda_4 
- \bar{\alpha}_2^2 \bar{\alpha}_3 \lambda_5 + \bar{\alpha}_2^3 \lambda_6 - \bar{\alpha}_0 \bar{\alpha}_2^2 \lambda_8,\\
&c_{9} = -\bar{\alpha}_3^3 - 6 \bar{\alpha}_0 \bar{\alpha}_3 \bar{\alpha}_6
+ (\bar{\alpha}_2 \bar{\alpha}_3^2 + 2 \bar{\alpha}_0 \bar{\alpha}_2 \bar{\alpha}_6) \lambda_1 
- \bar{\alpha}_2^2 \bar{\alpha}_3 \lambda_2 + \bar{\alpha}_2^3 \lambda_3 \\
&\quad \quad + 2 \bar{\alpha}_0 \bar{\alpha}_2 \bar{\alpha}_3 \lambda_4 - \bar{\alpha}_0 \bar{\alpha}_2^2  \lambda_5. 
\end{align*}
\end{enumerate}
Assign $\mathcal{H}$ to  $\widetilde{\mathcal{H}}$,
and  $\mathcal{I}$ to $\mathcal{I}^\ast$, then
return to the Step 1.

\medskip
Now we illustrate the algorithm with

\begin{exam}\label{E:RedAlgC34}
On the curve \eqref{C34Ex1} we implement the standard cryptography problem: add a degree $3$ divisor $D_3$ to itself 
and repeat this process iteratively. Two setups are considered:
(a) $D_3$  consists of three distinct points: $D_3= (-1,\,5)+(3,\,1)+(4,\,-3)$, and 
(b) $D_3$ consists of three equal points: $D_3=3(2,\,-1)$.
In each case we start the first iteration from Step 2,  
constructing $\mathcal{F}$ and $\mathcal{G}$ defining a divisor $D_4$
directly from the points.

\subsection*{Case (a): $n\, \{(-1,\,5)+(3,\,1)+(4,-3)\}$.}
Let $D_4 = 2(-1,\,5)+(3,\,1)+(4,-3)$, then $\mathcal{G}(x,y)$ is
defined by \eqref{GSame2Dist2P}, and  $\mathcal{F}(x) = (x+1)^2 (x-3)(x-4)$,
see Example~\ref{E:I2Same2DistPC1}. 
Then points $(3,\,1)$
and $(4,\,-3)$ are added consecutively, and so $2D_3$ divisor is obtained.
\begin{align*}
&D_3& 
&\mathcal{H}(x) = x^3 - 6 x^2 + 5 x + 12, \\
&&
&\mathcal{I}(x,y) = 3 x^2 + 5 y - x - 29; \\
&2D_3& 
&\mathcal{H}(x) = 26\,508\,107\,741\,423\,463\, x^3 + 397\,689\,267\,124\,379\,826\, x^2  \\
&& &\qquad\quad + 1\,870\,379\,369\,356\,664\,330\, x + 1\,524\,957\,157\,645\,063\,832, \\
&&
&\mathcal{I}(x,y) = 29\,634\,519\,963\, x^2 + 21\,261\,894\,998\, y + 289\,874\,841\,784\, x  \\
&& &\qquad\quad  + 165\,783\,154\,578; \\
&3D_3& 
&\mathcal{H}(x) =  132\,854\,002\,684\,150\,020\,787\,377\,323\,493\,369\,010\,137\,224\, x^3  \\
&& &\qquad\quad  -110\,131\,151\,863\,240\,535\,734\,181\,632\,313\,620\,530\,154\,204\,383\, x^2  \\
&& &\qquad\quad +  4\,240\,086\,218\,537\,526\,261\,755\,671\,705\,038\,721\,403\,541\,074\,994\, x \\
&& &\qquad\quad  -11\,846\,346\,070\,578\,926\,743\,102\,916\,547\,271\,037\,839\,019\,230\,108, \\
&&
&\mathcal{I}(x,y) =  2\,603\,652\,953\,486\,947\,241\,363\,417\,476\, x^2 \\
&& &\qquad\quad  + 232\,317\,982\,116\,414\,482\,349\,653\,949\,231\, y \\
&& &\qquad\quad  + 253\,588\,066\,753\,353\,560\,448\,633\,339\,503\, x \\
&& &\qquad\quad  -1\,766\,046\,605\,251\,168\,704\,501\,340\,315\,799.
\end{align*}

\subsection*{Case (b): $n \{3 (2,\,-1)\}$.} Let  $D_4 = 4(2,-1)$,
then $\mathcal{G}(x,y)$ is defined by \eqref{GSame4P}, and $\mathcal{F}(x) = (x-2)^4$,
see Example~\ref{E:I4SamePC1}. Then we add the point $(2,-1)$ consecutively. 
\begin{align*}
&D_3 & 
&\mathcal{H}(x) = x^3 - 6 x^2 + 12 x - 8, \\
&&
&\mathcal{I}(x,y) =  2\,162\, x^2 + 243\, y  - 10\,457\, x + 12\,509; \\
&2D_3 & 
&\mathcal{H}(x) = 1\,535\,080\,134\,464\,797\,694\,295\,656\,567\,778\,282\,346\,829\,227\,823\,678\,369\,398\,784 x^3 \\
&& &\qquad\quad -11\,654\,450\,685\,673\,839\,326\,280\,964\,785\,946\,543\,168\,248\,358\,809\,818\,717\,238\,849 x^2\\
&& &\qquad\quad + 6\,928\,653\,622\,424\,446\,486\,701\,046\,338\,242\,867\,004\,910\,441\,486\,486\,009\,884\,432 x \\ 
&& &\qquad\quad -1\,374\,416\,320\,134\,567\,210\,929\,966\,212\,279\,180\,215\,193\,999\,855\,041\,857\,833\,006, \\
&&
&\mathcal{I}(x,y) = 133\,072\,209\,440\,609\,507\,290\,563\,819\,791\,953\,601\,536 x^2 \\
&& &\qquad\quad  - 23\,240\,322\,439\,350\,339\,511\,372\,089\,021\,339\,773\,151 y \\
&& &\qquad\quad  -  896\,022\,294\,764\,013\,888\,618\,903\,959\,928\,856\,245\,643 x \\
&& &\qquad\quad + 190\,762\,431\,483\,770\,076\,812\,777\,831\,493\,327\,372\,463.
\end{align*}
\end{exam}

\section{Conclusion and discussion}
The proposed algorithm has essential advantages. It is direct and explicit, that is 
all computations within the reduction algorithm can be done in symbolic form
before the implementation. On curves $C_{3,4}$ and $C_{3,5}$ these computations
are presented in the paper, and the cases of curves $C_{3,7}$, $C_{3,8}$, $C_{3,10}$, $C_{3,11}$, 
$C_{3,13}$, $C_{3,14}$
are also completed by the authors. 
The algorithm has no upper bound of the number of points to add.

\appendix

\section{Reduction of a degree $g+1$ divisor on $C_{3,3\mFr+2}$}\label{A:Main35Family}
Here we give a solution of the reduction problem on a curve of the type~$\mathfrak{C}^2$.
It is similar to what is presented in Section~\ref{s:Main}. Below we suggest updates 
of subsections~\ref{ss:gg1Red} and \ref{ss:g1gRed}, whereas subsections
\ref{ss:g1Div} and \ref{ss:gInv} does not require any update.

\subsection{Minimal function defining a degree $g+1$ divisor}\label{A:35FamilyG}
Let a degree $g+1$ divisor be $D_{g+1} = D_g + (\bar{x},\, \bar{y})$, where a degree $g$ 
strictly non-special divisor $D_g$ is defined by $\mathcal{H}$ and  $\mathcal{I}$, according to Theorem~\ref{T:DgIG},
and $(\bar{x},\, \bar{y})$ does not coincide with any point from $D_g$.
Then $D_{g+1}$ is defined by
\begin{gather*}
\mathcal{F}(x) = \mathcal{H}(x) (x-\bar{x}),
\end{gather*} 
and an entire function $\mathcal{G}$ of weight $2g+1$, according to Theorem~\ref{T:g1FGCns}.

In order to find $\mathcal{G}$, we use \eqref{GfromI}, where  $\mathcal{I}$ has the form  \eqref{ImonomDef2}
and
\begin{align*}
& \mathcal{H}(x) = \sum_{k=0}^{3\mFr} h_{9\mFr-3k} x^k,\\
& \mathcal{M}(x) =  \sum_{i=0}^{3\mFr} a_{9\mFr-3i} x^i,\\
& \mathcal{N}(x,y) = y \sum_{i=0}^{\mFr-1} b_{3\mFr-3-3i} x^i + \sum_{i=0}^{2\mFr-1} b_{6\mFr-1-3i} x^i
\equiv y \nu_y(x) + \nu_x(x).
\end{align*}

Equating to zero the coefficients of the monomials $\{y x^{4\mFr}$, \ldots, $yx^{\mFr+1}\}$ and 
$\{x^{5\mFr}$, \ldots, $x^{2\mFr+2}\}$, we obtain $6\mFr-1$ linear equations of the form \eqref{BpEqs}
in $6\mFr + 1$ unknown parameters 
$$\mathfrak{p} = \big(a_{0},\, a_3,\, \dots,\, a_{9\mFr},\, b_0,\, b_3,\, \ldots,\, b_{3\mFr - 3},\, 
b_2,\, b_5,\, \ldots,\, b_{6\mFr - 1}\big).$$ 
The coefficient matrix has the block form \eqref{BlockM}.
Block $\mathfrak{G}_a$ has size $3\mFr \times (3\mFr +1)$ with one zero column on the right,
the remaining part is  lower triangular $(\mFr +1)$-diagonal, namely:
\begin{gather*}
 \mathfrak{G}_a = 
 \begin{pmatrix} \alpha_0 & 0 & \cdots & 0 & 0 & \cdots & 0 & 0 \\
\alpha_3 &  \alpha_0 & \ddots & \vdots & \vdots &  & \vdots & \vdots  \\
\vdots & \ddots & \ddots & 0 & 0 & \cdots & 0 & 0  \\
\alpha_{3 \mFr}  &  \ddots & \alpha_3 &  \alpha_0 & 0 &  \cdots & 0 & 0 \\
0 & \alpha_{3 \mFr}  & \ddots & \alpha_3 & \ddots & \ddots & \vdots & \vdots \\
\vdots & & \ddots & \ddots & \ddots &   \alpha_0 & 0 & 0  \\
0 & \cdots & 0 & \alpha_{3 \mFr} &  \cdots & \alpha_3 &  \alpha_0 & 0 
 \end{pmatrix}.
\end{gather*}
Block $\mathfrak{G}_b$ is of size $(3\mFr - 1) \times (3\mFr + 1)$ with two zero columns on the right,
the remaining part is lower triangular $(2 \mFr + 1)$-diagonal:
\begin{gather*}
\mathfrak{G}_b = 
\begin{pmatrix} \alpha_2 & 0 & \cdots & 0 & 0 & \cdots & 0 & 0 & 0\\
\alpha_5 &  \alpha_2 & \ddots & \vdots & \vdots &  & \vdots & \vdots & \vdots \\
\vdots & \ddots & \ddots & 0 & 0 & \cdots & 0 & 0 & 0\\
\alpha_{6\mFr + 2}  & \ddots & \alpha_5 &  \alpha_2 & 0  & \cdots & 0 & 0 & 0 \\
0 & \alpha_{6\mFr + 2} & \ddots & \alpha_5 & \ddots &  \ddots & \vdots & \vdots & \vdots \\
\vdots & & \ddots & \ddots &\ddots &\alpha_2 &0 & 0 & 0 \\
0  & \cdots & 0 & \alpha_{6\mFr + 2}  & \cdots & \alpha_5 &  \alpha_2 & 0 & 0
 \end{pmatrix}.
\end{gather*}
Block $\mathfrak{H}_a$ is square of size $3\mFr$ with $2\mFr$ zero columns on the right, 
and block $\mathfrak{H}_b$ is of size $(3\mFr-1) \times 3\mFr$ with $\mFr$ zero columns on the left
\begin{gather*}
\mathfrak{H}_a =  \begin{pmatrix} 
h_0 & 0 & \cdots & 0 & 0 & \cdots & 0 \\
h_3 & h_0 & \ddots & \vdots & 0 & \cdots & 0 \\
\vdots & \ddots & \ddots & 0 & \vdots & & \vdots \\
h_{3\mFr-3} & h_{3\mFr - 6} & \ddots & h_0 & 0 & \cdots & 0 \\
h_{3\mFr} & h_{3\mFr - 3}  & \cdots & h_3 & 0 &  \cdots & 0 \\
\vdots & \vdots & & \vdots & \vdots & & \vdots \\
h_{9\mFr - 3} & h_{9\mFr - 6} & \cdots & h_{6\mFr} & 0 & \cdots &0\end{pmatrix},\\
\mathfrak{H}_b =  \begin{pmatrix} 
0 & \cdots & 0 & h_0 & 0 & \cdots & 0 \\
0 & \cdots & 0 & h_3 & h_0 & \cdots & 0\\
\vdots & &  \vdots & \vdots & \vdots & \ddots & \vdots \\
0 & \cdots & 0 & h_{6\mFr - 3} & h_{6\mFr - 6} & \cdots & h_0\\
\vdots & &  \vdots & \vdots & \vdots &  & \vdots \\
0 & \cdots & 0 & h_{9\mFr - 6} & h_{9\mFr - 9} & \cdots & h_{3\mFr - 3}
\end{pmatrix}.
\end{gather*}
The kernel of $\mathfrak{B}$ is two-dimensional; and  $a_{9\mFr}$ is a free parameter, since
the $(3\mFr+1)$-th column is zero. Also $\mathfrak{p}$ is defined up to a constant multiple, which is $a_0$.
We introduce $\widetilde{\mathcal{G}}(x,y)$
with $a_{9\mFr}=0$.  Then assign
\begin{gather*}
a_{9\mFr} = -\frac{\widetilde{\mathcal{G}}(\bar{x},\bar{y})}{\mathcal{I}(\bar{x},\bar{y})}.
\end{gather*}
Finally, $\mathcal{G}$ is constructed by \eqref{Gexpr}.

\subsection{Minimal function defining a degree $g$ divisor}\label{A:35FamilyI}
Let $D^\ast_g$ be a degree $g$ strictly non-special divisor  
defined by a polynomial $\mathcal{H}^\ast$ of degree $g$ and an entire
function $\mathcal{G}$ of weight $2g+1$.  We
find the minimal function $\mathcal{I}$ defining $D^\ast_g$. 
For this purpose we use \eqref{IfromG},
where $\mathcal{G}$ has the form \eqref{GmonomDef2}, and
\begin{align*}
& \mathcal{H}^\ast (x) = \sum_{k=0}^{3\mFr+1} \bar{h}_{9\mFr+3-3k} x^k,\\
& \mathcal{M}(x) =  \sum_{i=0}^{3\mFr} a_{9\mFr-3i} x^i,\\
& \mathcal{N}(x,y) = y \sum_{i=0}^{\mFr-1} b_{3\mFr-2-3i} x^i + \sum_{i=0}^{2\mFr} b_{6\mFr-3i} x^i
\equiv y \nu_y(x) + \nu_x(x).
\end{align*}

Equating to zero the coefficients of the monomials $\{y x^{4\mFr}$, \ldots, $yx^{\mFr+1}\}$ and 
$\{x^{5\mFr+1}$, \ldots, $x^{2\mFr+1}\}$, we obtain $6\mFr+1$ linear equations of the form \eqref{BpEqs}
in $6\mFr + 2$ unknown parameters 
$$\mathfrak{p} = \big(a_{0},\, a_3,\, \dots,\, a_{9\mFr},\, b_{1},\, b_4,\, \dots,\, b_{3\mFr - 2},\, 
b_{0},\, b_3,\, \dots,\, b_{6\mFr} \big).$$ 
The coefficient matrix has the block form \eqref{BlockM}.
Block $\mathfrak{G}_a$ has size $3\mFr \times (3\mFr +1)$ with one zero column on the right,
the remaining part is  lower triangular $(\mFr +1)$-diagonal, namely:
\begin{gather*}
\mathfrak{G}_a = 
 \begin{pmatrix} \gamma_1 & 0 & \dots & 0 & 0 & \dots & 0 & 0 \\
\gamma_4 &  \gamma_1 & \ddots & \vdots & \vdots &  & \vdots & \vdots \\
\vdots & \ddots & \ddots & 0 & 0 & \dots & 0 & 0 \\
\gamma_{3\mFr + 1}  &  \ddots & \gamma_4 &  \gamma_1 & 0 & \dots &  0& 0\\
0 & \gamma_{3\mFr + 1} & \ddots & \gamma_4 & \ddots & \ddots & \vdots & 0 \\
\vdots & & \ddots & \ddots & \ddots &  \gamma_1 & 0 & \vdots \\
0 & \dots & 0 & \gamma_{3\mFr + 1} &  \dots & \gamma_4 &  \gamma_1 &0
 \end{pmatrix}.
\end{gather*}
Block $\mathfrak{G}_b$ is square of size $(3\mFr+1)$, and lower triangular $(2\mFr + 2)$-diagonal
\begin{gather*}
\mathfrak{G}_b = 
\begin{pmatrix} \gamma_0 & 0 & \dots & 0 & 0 & \dots & 0  \\
\gamma_3 &  \gamma_0 & \ddots & \vdots & \vdots & &  \vdots   \\
\vdots & \ddots & \ddots & 0 & 0 & \dots & 0  \\
\gamma_{6 \mFr+3}  &  \ddots & \gamma_3 &  \gamma_0 & 0 & \dots & 0 \\
0 & \gamma_{6 \mFr+3} & \ddots & \gamma_3 & \ddots & \ddots & \vdots  \\
\vdots & & \ddots & \ddots & \ddots &   \gamma_0  & 0 \\
0 & \dots & 0 & \gamma_{6 \mFr+3} &  \dots & \gamma_3 &  \gamma_0 
 \end{pmatrix}.
\end{gather*}
Block $\mathfrak{H}_a$ of size $3\mFr \times (3\mFr + 1)$ has $2\mFr + 1$ zero columns on the right,
and block  $\mathfrak{H}_b$ of size $(3\mFr + 1)\times (3 \mFr + 1)$ has $\mFr$ zero columns on the left, namely:
\begin{gather*}
\mathfrak{H}_a =  \begin{pmatrix} 
h_0 & 0 & \dots & 0 & 0 & & 0 \\
h_3 & h_0 & \ddots & \vdots & 0 & \dots & 0 \\
\vdots & \ddots & \ddots & 0 & \vdots & & \vdots \\
h_{3\mFr - 3} & h_{3\mFr - 6} & \ddots & h_0 & 0 & \dots & 0 \\
h_{3\mFr} & h_{3\mFr - 3}  & \dots & h_3 & 0 &  \dots & 0 \\
\vdots & \vdots & & \vdots & \vdots & & \vdots \\
h_{9\mFr - 3} & h_{9\mFr - 6} & \dots & h_{6\mFr} & 0 & \dots &0\end{pmatrix},\\
\mathfrak{H}_b =  \begin{pmatrix} 
0 & \dots & 0 & h_0 & 0 & \dots & 0 \\
0 & \dots & 0 & h_3 & h_0 & \dots & 0\\
\vdots & &  \vdots & \vdots & \vdots & \ddots & \vdots \\
0 & \dots & 0 & h_{6\mFr} & h_{6\mFr - 3} & \dots & h_0\\
\vdots & &  \vdots & \vdots & \vdots & \ddots & \vdots \\
0 & \dots & 0 & h_{9 \mFr} & h_{9 \mFr - 3} & \dots & h_{6 \mFr}
\end{pmatrix}.
\end{gather*}
The kernel of $\mathfrak{B}$ is one-dimensional, so the parameters $\mathfrak{p}$
are defined up to a constant multiple, which is $a_0$.

\section{Implementation of reduction algorithm on $C_{3,5}$}\label{A:RedAlgC35}
Let a curve $C_{3,5}$ be defined by \eqref{C35eq1}.
We start each iteration, except possibly the first, with
\begin{align*}
&\mathcal{H}(x) = h_0 x^4 +  h_3 x^3 + h_6 x^2 + h_9 x + h_{12},\\
&\mathcal{I}(x,y) = \alpha_0 x y + \alpha_2 x^2 + \alpha_3 y  +  \alpha_5 x + \alpha_8,
\end{align*}
which define the reduce divisor $D_4$ from the previous iteration.
We denote the point which is added by $(\bar{x},\, \bar{y})$.
\begin{enumerate}
\item[\textbf{Step 1.}] Find
\begin{align*}
&\mathcal{F}(x) = f_0 x^5 + f_3 x^4 +  f_6 x^3 + f_9 x^2 + f_{12} x + f_{15},\\
&\widetilde{\mathcal{G}}(x,y) = \bar{\gamma}_0 x^3 + \bar{\gamma}_1 x y + \bar{\gamma}_3 x^2 
+ \bar{\gamma}_4 y  +  \bar{\gamma}_6 x + \bar{\gamma}_9,
\end{align*}
where 
\begin{equation*}
f_{0} = h_0,\qquad f_{3i+3} = h_{3i+3} - h_{3i} \bar{x},\ i=0,\,\dots,\,3, \qquad
f_{15} = - h_{12} \bar{x},
\end{equation*}
and
\begin{align*}
&\bar{\gamma}_0 = h_0,\\
&\bar{\gamma}_1 = \sfit{C}_8^{-1} \big(h_0 \alpha_3^3 - h_3 \alpha_0 \alpha_3^2 
+ h_6 \alpha_0^2 \alpha_3 - h_9 \alpha_0^3 \big),\\
&\bar{\gamma}_3 = \sfit{C}_8^{-1} \big( h_0 \alpha_3 (\alpha_3 \alpha_5 - \alpha_0 \alpha_8)
- h_3 \alpha_0 (\alpha_3 \alpha_5 - \alpha_0 \alpha_8) 
 + h_6 \alpha_0 \alpha_2 \alpha_3
- h_9 \alpha_0^2 \alpha_2 \big), \\
&\bar{\gamma}_4 =  - \sfit{C}_8^{-1}  h_{12} \alpha_0^3,\\
&\bar{\gamma}_6 = \sfit{C}_8^{-1} \big( h_0 \alpha_3^2 \alpha_8
- h_3 \alpha_0 \alpha_3 \alpha_8  + h_6 \alpha_0^2 \alpha_8
- h_9 \alpha_0 (\alpha_0 \alpha_5 - \alpha_2 \alpha_3)  - h_{12} \alpha_0^2 \alpha_2 \big),\\
&\bar{\gamma}_9 =  \sfit{C}_8^{-1}  h_{12} \alpha_0 (\alpha_2 \alpha_3 - \alpha_0 \alpha_5), \\
&\sfit{C}_8 = \alpha_0^2  \alpha_8 - \alpha_0 \alpha_3 \alpha_5 + \alpha_2 \alpha_3^2.
\end{align*}
Let $c = \widetilde{\mathcal{G}}(\bar{x},\bar{y})/ \mathcal{I}(\bar{x},\bar{y})$. Then
\begin{align*}
&\mathcal{G}(x,y) = \gamma_0 x^3 + \gamma_1 x y + \gamma_3 x^2 
+ \gamma_4 y  +  \gamma_6 x + \gamma_9,
\end{align*}
where
\begin{align*}
&\gamma_0 = \bar{\gamma}_0, \qquad
\gamma_{n} = \bar{\gamma}_{n} - c \alpha_{n-1},\quad n = 1,\, 3,\, 4,\, 6,\, 9.
\end{align*}

\item[\textbf{Step 2.}] 
Find
\begin{equation*}
\mathcal{H}^\ast (x) = \bar{h}_0 x^4 + \bar{h}_{3} x^3 + \bar{h}_{6} x^2 + \bar{h}_{9} x + \bar{h}_{12},
\end{equation*}
with
\begin{gather*}
\bar{h}_{12-3 k} = \sum_{i=0}^{k} (-1)^{i} d_{-15-3i} c_{27-3(k-i)},\qquad k=0,\, \dots,\, 4,
\end{gather*}
where 
\begin{gather*}
\begin{split} \notag
&d_{-15} = f_{15}^{-1},\\
&d_{-15-3j} = \sum_{i=1}^{j} (-1)^{i-1} f_{15}^{-1} f_{15-3i} d_{-15-3(j-i)},\qquad j = 1,\dots, 5,
\end{split}
\end{gather*}
and $c_{15}$, $c_{18}$, $c_{21}$, $c_{24}$, $c_{27}$ are coefficients of $\mathcal{Z}_{9}$ defined by \eqref{Zgp}:
\begin{align*}
&c_{27} = - \gamma_9^3 + \gamma_4 \gamma_9^2 \lambda_5  - \gamma_4^2 \gamma_9  \lambda_{10} 
+ \gamma_4^3 \lambda_{15},\\
&c_{24} =  - 3 \gamma_6 \gamma_9^2 + \gamma_4 \gamma_9^2 \lambda_2   
+ (2\gamma_4 \gamma_6 \gamma_9 + \gamma_1 \gamma_9^2) \lambda_5 
-  \gamma_4^2 \gamma_9  \lambda_7 \\
&\quad \quad
- (\gamma_4^2 \gamma_6 + 2 \gamma_1 \gamma_4 \gamma_9) \lambda_{10}
+ \gamma_4^3 \lambda_{12} + 3 \gamma_1 \gamma_4^2 \lambda_{15},\\
&c_{21} = -3 \gamma_6^2 \gamma_9 - 3 \gamma_3 \gamma_9^2
+ (2\gamma_4 \gamma_6 \gamma_9 + \gamma_1 \gamma_9^2) \lambda_2 
- \gamma_4^2 \gamma_9 \lambda_4  \\
&\quad \quad
+ (\gamma_4 \gamma_6^2 + 2 \gamma_3 \gamma_4 \gamma_9 
+ 2 \gamma_1 \gamma_6 \gamma_9) \lambda_5
- (\gamma_4^2 \gamma_6 + 2 \gamma_1 \gamma_4 \gamma_9) \lambda_7   + \gamma_4^3 \lambda_9\\
&\quad \quad
- (\gamma_3 \gamma_4^2 + 2 \gamma_1 \gamma_4 \gamma_6 
+ \gamma_1^2 \gamma_9) \lambda_{10} + 3 \gamma_1 \gamma_4^2  \lambda_{12} 
+ 3 \gamma_1 \gamma_4 \lambda_{15},\\
&c_{18} =  - \gamma_6^3  - 6 \gamma_3 \gamma_6 \gamma_9 - 3 \gamma_0 \gamma_9^2
- \gamma_4^2 \gamma_9 \lambda_1  
 + (\gamma_4 \gamma_6^2 + 2 \gamma_3 \gamma_4 \gamma_9
 + 2 \gamma_1 \gamma_6 \gamma_9) \lambda_2 \\
&\quad \quad + 
(\gamma_4^2 \gamma_6+ 2 \gamma_1 \gamma_4 \gamma_9 ) \lambda_4 
+ (2\gamma_3 \gamma_4 \gamma_6 + \gamma_1 \gamma_6^2 
+ 2 \gamma_1 \gamma_3 \gamma_9 + 2 \gamma_0 \gamma_4 \gamma_9) \lambda_5\\
&\quad \quad 
- (\gamma_3 \gamma_4^2 + 2 \gamma_1 \gamma_4 \gamma_6 
+ \gamma_1^2\gamma_9) \lambda_7 + 3 \gamma_1 \gamma_4^2 \lambda_9\\
&\quad\quad  - (2\gamma_1\gamma_3 \gamma_4 +\gamma_0 \gamma_4^2 
+ \gamma_1^2 \gamma_6) \lambda_{10} 
+ 3 \gamma_1 \gamma_4 \lambda_{12} + \gamma_1^3\lambda_{15},\\
&c_{15} = -3 \big(\gamma _9 \gamma
   _3^2+\gamma _6^2 \gamma _3+2 \gamma _0 \gamma _6 \gamma _9\big) -
   \big(\gamma _4^2 \gamma _6+2 \gamma _1 \gamma _4 \gamma_9\big) \lambda _1
   +\big(\gamma _1 \gamma
   _6^2+2 \gamma _3 \gamma _4 \gamma _6 \\
   &\quad\quad  +2 \gamma _1
   \gamma _3 \gamma _9+2 \gamma _0 \gamma _4 \gamma
   _9\big) \lambda _2 + \gamma _4^3 \lambda_3 
    - \big(\gamma _9 \gamma _1^2 + 2 \gamma _4
   \gamma _6 \gamma _1 + \gamma _3 \gamma _4^2\big) \lambda _4\\
   &\quad\quad
   + \big(\gamma _4 \gamma _3^2+2 \gamma
   _1 \gamma _6 \gamma _3+2 \gamma _0 \gamma _4
   \gamma _6+2 \gamma _0 \gamma _1 \gamma _9\big) \lambda _5  
   - \big(\gamma _6 \gamma _1^2 + 2 \gamma_3 \gamma _4 \gamma _1 
   + \gamma _0 \gamma_4^2\big) \lambda _7 \\
   &\quad\quad 
    + \gamma _1^3 \lambda _{12} + 3 \gamma _4 \gamma _1^2  \lambda _9 
  + 3 \gamma _4^2 \gamma _1 \lambda
   _6-\big(\gamma _1 \gamma _3+2 \gamma _0 \gamma
   _4\big) \gamma _1 \lambda _{10}.
\end{align*}

\item[\textbf{Step 3.}] Find the new
\begin{equation*}
\mathcal{I}^\ast (x,y) = \bar{\alpha}_0 x y + \bar{\alpha}_2 x^2 + \bar{\alpha}_3 y  
+  \bar{\alpha}_5 x + \bar{\alpha}_8,
\end{equation*}
where
\begin{align*}
&\bar{\alpha}_0 = \bar{h}_0,\\
&\bar{\alpha}_2 = \sfit{C}_{12}^{-1}\big(-\gamma
   _1^2 \gamma _0^2 \bar{h}_{12} +\gamma _1 \gamma _4 \gamma _0^2 \bar{h}_9
  - \gamma _4^2 \gamma _0^2 \bar{h}_6+ (\gamma _0 \gamma _9 \gamma _1^2-\gamma
   _0 \gamma _4 \gamma _6 \gamma _1+\gamma _0 \gamma
   _3 \gamma _4^2) \bar{h}_3 \\
   &\quad\quad +(\gamma _3 \gamma _4 \gamma _6
   \gamma _1 - \gamma _3 \gamma
   _9 \gamma _1^2-\gamma _0 \gamma _4 \gamma _9 \gamma
   _1-\gamma _3^2 \gamma _4^2+\gamma _0 \gamma _4^2
   \gamma _6) \bar{h}_0 \big),\\
&\bar{\alpha}_3 = \sfit{C}_{12}^{-1}\big(-\gamma _0
   \gamma _1^3 \bar{h}_{12} + \gamma _0 \gamma
   _4 \gamma _1^2 \bar{h}_9 - \gamma _0 \gamma _4^2 \gamma
   _1 \bar{h}_6 + \gamma _0 \gamma _4^3 \bar{h}_3 + (\gamma _1 \gamma _6 \gamma
   _4^2 - \gamma _3 \gamma _4^3 - \gamma _1^2 \gamma _9 \gamma _4) \bar{h}_0 \big), \\
&\bar{\alpha}_5 = \sfit{C}_{12}^{-1}\big((\gamma _0^2 \gamma _1
   \gamma _4-\gamma _0 \gamma _1^2 \gamma _3)
   \bar{h}_{12}+(\gamma _0 \gamma _1 \gamma _3 \gamma
   _4-\gamma _0^2 \gamma _4^2)
   \bar{h}_9+(\gamma _0 \gamma _1^2 \gamma _9-\gamma
   _0 \gamma _1 \gamma _4 \gamma _6)
   \bar{h}_6 \\ &\quad\quad +(\gamma _0 \gamma _4^2 \gamma _6-\gamma
   _0 \gamma _1 \gamma _4 \gamma _9)
   \bar{h}_3+(\gamma
   _4 \gamma _6^2 \gamma _1 - \gamma _6 \gamma _9 \gamma _1^2 -
    \gamma _3 \gamma _4^2 \gamma _6+\gamma _0 \gamma _4^2 \gamma _9)
   \bar{h}_0 \big),\\
&\bar{\alpha}_8 = \sfit{C}_{12}^{-1}\big((\gamma_0 \gamma _3 \gamma _4 \gamma _1 
   -\gamma _0 \gamma _6 \gamma _1^2-\gamma _0^2
   \gamma _4^2) \bar{h}_{12}
   + \gamma _0 \gamma _9 \gamma _1^2
   \bar{h}_9 - \gamma _0 \gamma _4 \gamma _9 \gamma _1
   \bar{h}_6 + \gamma _0 \gamma _4^2
   \gamma _9 \bar{h}_3 \\ &\quad\quad 
   + (\gamma _1 \gamma _6 \gamma _9 \gamma
   _4 -\gamma _3 \gamma _9 \gamma
   _4^2-\gamma _1^2 \gamma _9^2) \bar{h}_0\big), \\
&\sfit{C}_{12} = \gamma_0 \gamma_4^3 - \gamma_1 \gamma_3 \gamma_4^2
+ \gamma_1 \gamma_4 \gamma_6 - \gamma_1^3 \gamma_9. 
\end{align*}

\item[\textbf{Step 4.}] Find
\begin{equation*}
\widetilde{\mathcal{H}}^\ast (x) = \tilde{h}_0 x^4 + \tilde{h}_{3} x^3 + \tilde{h}_{6} x^2 + \tilde{h}_{9} x + \tilde{h}_{12},
\end{equation*}
with
\begin{gather*}
\tilde{h}_{12-3 k} = \sum_{i=0}^{k} (-1)^{i} d_{-12-3i} c_{24-3(k-i)},\qquad k=0,\, \dots,\, 4,
\end{gather*}
where 
\begin{gather*}
\begin{split} \notag
&d_{-12} = \bar{h}_{12}^{-1},\\
&d_{-12-3j} = \sum_{i=1}^{j} (-1)^{i-1} \bar{h}_{12}^{-1} \bar{h}_{12-3i} d_{-12-3(j-i)},\qquad j = 1,\dots, 4,
\end{split}
\end{gather*}
and $c_{12}$, $c_{15}$, $c_{18}$, $c_{21}$, $c_{24}$ are coefficients of $\mathcal{Z}_{8}$ defined by \eqref{Zgp}:
\begin{align*}
&c_{24} = \alpha _3^3\lambda _{15} - \alpha _8
   \alpha _3^2  \lambda _{10} + \alpha _8^2 \alpha _3  \lambda _5 - \alpha _8^3,\\
&c_{21} = 3 \alpha _0 \alpha _3^2  \lambda _{15} 
 + \alpha_3^3 \lambda _{12} 
  - (\alpha _3^2 \alpha _5+2 \alpha _0 \alpha _3 \alpha _8) \lambda _{10} 
  - \alpha _8 \alpha _3^2 \lambda _7  \\
   &\quad\quad 
   +\alpha _8 (2 \alpha _3 \alpha_5 + \alpha _0 \alpha _8) \lambda_5 + \alpha _8^2 \alpha _3 \lambda _2
   - 3 \alpha_5 \alpha _8^2,\\
&c_{18} = 3 \alpha_0^2 \alpha _3 \lambda _{15} + 3 \alpha _0 \alpha _3^2 \lambda _{12}
- (\alpha _8 \alpha_0^2 + 2 \alpha _3 \alpha _5 \alpha _0 + \alpha _2 \alpha _3^2) \lambda _{10} \\
  &\quad\quad
 + \alpha _3^3 \lambda _9 -(\alpha _3^2 \alpha _5+2 \alpha _0 \alpha _3 \alpha _8) \lambda _7  
  + (\alpha _3 \alpha _5^2+2 \alpha _0 \alpha_8 \alpha _5
  + 2 \alpha _2 \alpha _3 \alpha_8) \lambda _5 \\
 &\quad\quad  - \alpha _8 \alpha _3^2 \lambda _4 
   + \alpha _8 (2 \alpha _3 \alpha_5 + \alpha _0 \alpha _8) \lambda_2
    - 3 \alpha _8 (\alpha _5^2+\alpha _2 \alpha_8), \\
&c_{15} = \alpha _0^3  \lambda_{15}
 + 3  \alpha _0^2 \alpha _3 \lambda _{12}
   - (2 \alpha _0 \alpha _2 \alpha _3 + \alpha _0^2 \alpha_5) \lambda _{10}  
     +3 \alpha _3^2  \alpha_0 \lambda _9 \\
     &\quad\quad
     -(\alpha _8 \alpha _0^2 + 2 \alpha _3 \alpha _5 \alpha _0 + \alpha _2 \alpha _3^2)
   \lambda _7
  +2 \alpha _0 \alpha _2 \alpha_8) \lambda _5+\alpha _3^3 \lambda_6 \\
   &\quad\quad
     +(\alpha _0 \alpha _5^2+2 \alpha _2 \alpha_3 \alpha _5
    -\alpha _3 (\alpha _3
   \alpha _5+2 \alpha _0 \alpha _8) \lambda_4 \\
   &\quad\quad
   +(\alpha _3 \alpha _5^2+2 \alpha _0 \alpha
   _8 \alpha _5+2 \alpha _2 \alpha _3 \alpha_8) \lambda _2
       - \alpha _3^2 \alpha _8 \lambda_1    
      - (\alpha _5^3  +6 \alpha _2 \alpha _5  \alpha_8), \\
&c_{12} =  \alpha _0^3 \lambda _{12} 
   - \alpha _2 \alpha _0^2 \lambda _{10} 
    + 3 \alpha _3\alpha _0^2  \lambda _9  
    -(2 \alpha _0 \alpha _2 \alpha _3+\alpha_0^2 \alpha _5)  \lambda _7 
    + 3 \alpha _3^2  \alpha _0 \lambda_6 \\
    &\quad\quad
   + (\alpha _2^2 \alpha _3+2
   \alpha _0 \alpha _2 \alpha _5) \lambda _5
   -(\alpha _8 \alpha _0^2 + 2 \alpha _3
   \alpha _5 \alpha _0 + \alpha _2 \alpha _3^2)
   \lambda _4
    + \alpha _3^3 \lambda _3 \\
    &\quad\quad
    + (\alpha _0 \alpha _5^2+2 \alpha _2 \alpha
   _3 \alpha _5+2 \alpha _0 \alpha _2 \alpha
   _8) \lambda _2 
   - (\alpha _3^2 \alpha _5 + 2 \alpha _0 \alpha _3 \alpha _8) \lambda_1 
   - 3  (\alpha_2 \alpha _5^2+\alpha _2^2 \alpha_8).
\end{align*}

\end{enumerate}

\section{Division of polynomials}\label{A:PolyDiv}
Here a division of polynomial $\mathcal{Z}$ by polynomial $\mathcal{F}$ is given explicitly in terms of their coefficients.
For the purpose of the present paper the quotient polynomial $\mathcal{H}$ is supposed of degree $g$, and
coefficients are labeled by their S\={a}to weights. Denote
\begin{gather*}
\mathcal{Z} = \sum_{i=0}^N c_{3(N-i)} x^i = c_0 x^N + c_3 x^{N-1} + \cdots + c_{3N-3} x +c_{3N},\\
\mathcal{F} = \sum_{i=0}^n f_{3(n-i)} x^i = f_0 x^n + f_3 x^{n-1} + \cdots + f_{3n-3} x + f_{3n},\\
\mathcal{H} = \sum_{i=0}^g h_{3g-3i} x^i = h_0 x^g + h_3 x^{g-1} + \cdots + h_{3g-3} x + h_{3g},
\end{gather*} 
where $N=g+n$, and $\mathcal{Z}$ is divisible by $\mathcal{F}$.
Then
\begin{subequations}\label{PDivCoefs}
\begin{gather}
h_{3(g-k)} = \sum_{i=0}^{k} (-1)^{i} d_{-3(N-g+i)} c_{3(N-k+i)},
\end{gather}
and
\begin{gather}\label{dExpr}
d_{-3(n+j)} = \sum_{i=1}^{j} (-1)^{i-1} f_{3n}^{-1} f_{3(n-i)} d_{-3(n+j-i)}.
\end{gather}
\end{subequations}


\end{document}